\documentstyle{amsppt}

\topmatter

\author Richard Haydon \endauthor
\title Trees in renorming theory \endtitle
\address Brasenose College, Oxford, OX1 4AJ \endaddress
\email richard.haydon\@brasenose.oxford.ac.uk \endemail

\endtopmatter

\document

\define\Beta{\text{\rm B}}
\define\Alpha{\text{\rm A}}

\define\C{{\Cal C}}

\define\N{{\Cal N}}

\define\reals{{\Bbb R}}

\define\rats{{\Bbb Q}}
\define\nats{{\Bbb N}}

\define\dual{{^{\textstyle*}}}
\define\norm#1{{\|#1\|}}
\define\onorm#1{{\|#1\|_{\text{ord}}}}
\define\anorm#1{{\|#1\|_{\text{arb}}}}
\define\inorm#1{{\|#1\|_{\infty}}}
\define\half{{\textstyle\frac12}}

\define\quarter{{\textstyle\frac14}}

\define\sixth{{\textstyle\frac16}}

\define\supp{\text{\rm supp}\,}

\define\osc{\text{\rm osc}\,}

\define\ball{\text{\rm ball}\,}
\define\dom{\text{\rm dom}\,}
\define\rg{\text{\rm im}\,}
\define\less{\setminus}

\define\indic{{\pmb 1}}
\define\rest#1{\restriction_{\textstyle #1}}

\define\<{\langle}
\define\>{\rangle}

\define\trg{\succ}
\define\trge{\succcurlyeq}
\define\trl{\prec}
\define\trle{\preccurlyeq}

\define\tp{\tau _{\text{p}}}

\define\Choquet{1}
\define\Day {2}
\define\DGZ{3}
\define\DGZtwo{4}
\define\EkLeb{5}
\define\Fremlin{6}
\define\GTWZ{7}
\define\Haydon{8}
\define\Haydonscat{9}
\define\Haydontwo{10}
\define\Haydonthree{11}
\define\Haydonfour{12}
\define\Haydonfive{13}
\define\JNR{14}
\define\NamPol{15}
\define\Preiss{16}
\define\Smulyan{17}
\define\Tal{18}
\define\Tod{19}
\define\Troy{20}

\heading Introduction \endheading

Renorming theory is that branch of functional analysis which investigates 
problems of the form: for which Banach spaces $X$ does there exist a norm on 
$X$, equivalent to the given norm, with some good geometrical 
property of smoothness or strict convexity? The hope is to give answers in 
terms of familiar linear topological properties such as reflexivity, 
separability, separability of the dual space $X\dual$ and so on. The 
philosophy is well summed up by the title of Section VII.4 of Day's book 
[\Day], ``Isomorphisms to improve the norm''. An up-to-date account of the 
theory is given in the recent and authoritative text of Deville, Godefroy and 
Zizler [\DGZ]. 

In the next section we shall give full definitions of all the 
smoothness and convexity conditions to be considered. For the moment, let us 
recall that a norm $\norm{\cdot}$ on a Banach space $X$ is said to be 
Fr\'echet-smooth (F) if the function $x\mapsto \norm x$ is 
Fr\'echet-differentiable except at 0. We say that $\norm\cdot$ is locally 
uniformly convex (LUR) if, whenever $x$ and $x_n$ $(n\in \nats)$ are elements 
of $X$ such that $\norm {x_n}\to \norm x$ and $\norm{x+x_n}\to 2\norm x$ as 
$n\to \infty$, we necessarily have $\norm {x-x_n}\to 0$. One reason for being 
interested in norms with this property is that if the dual norm 
$\norm{\cdot}\dual$ on $X\dual$ is LUR then the original norm $\norm{\cdot}$ 
on $X$ is (F). 

A class of Banach spaces which has turned out to be of special importance is 
the class of strong differentiability spaces, now called {\sl Asplund spaces}. 
A Banach space $X$ is said to be an Asplund space if every continuous convex 
real-valued function, defined on a convex open subset of $X$, is 
Fr\'echet-differentiable at all points of a dense $\Cal G_\delta $ subset of 
its domain. Work of Stegall, Namioka and Phelps established remarkable 
equivalent characterizations of this class of space: $X$ is Asplund if and 
only if the dual space $X\dual$ has the Radon--Nikodym Property and if and 
only if every separable subspace of $X$ has separable dual. In particular, a
space $\C_0(L)$, with $L$ locally compact, is Asplund if and only if $L$ is 
{\sl scattered}. Preiss [\Preiss] recently established dense
differentiability for Lipschitz functions on Asplund spaces, and Asplund 
spaces have also proved to be of interest in Optimization Theory. 

In the case of a separable space $X$, there is a very tight connection between
renorming and the Asplund Property. Indeed, in this special case the following 
are all equivalent: separability of the dual space $X\dual$; the existence on
$X$ of an equivalent norm with LUR dual norm; the existence on $X$ of an 
equivalent norm which is (F). The question whether there are implications 
between Fr\'echet renormability and the Asplund Property for general Banach 
spaces was posed in [\Day]. 
One implication was established by Ekeland and Lebourg [\EkLeb] in the course 
of their work on perturbed optimization problems: if $X$ has an equivalent 
Fr\'echet-differentiable norm then $X$ is Asplund. That the converse does not 
hold was the main result of the author's paper [\Haydon].

The example presented in [\Haydon] was a Banach space $\C_0(L)$ with $L$ a
locally compact scattered space of a special kind, namely a tree. From the 
point of view of renorming theory this Banach space was very bad indeed, 
admitting neither a strictly convex renorming nor a 
{\sl G\^ateaux\/}-differentiable renorming. In order to achieve this degree of 
bad behaviour, the tree had to be chosen to be very large and the author 
suggested at the end of [\Haydon] that it might be interesting to study 
renormings of the spaces of continuous functions on trees with a more subtle 
structure. The present paper is the result of that study. The aim has been, 
wherever possible, to establish necessary and sufficient conditions, expressed 
in terms of the combinatorial structure of a tree $\Upsilon $, for the 
existence of equivalent good norms on the space $\C_0(\Upsilon )$. We are then  
able to give counterexamples to a number of open questions about renormings 
by considering suitably chosen trees. 

As has been found in Logic and General Topology (see for instance 
[\Tod]), trees are very agreeable objects to work with, offering a diversity 
of behaviour within a structure that is sufficiently simple to admit precise 
analysis. Thus we are able to offer fairly satisfactory necessary and 
sufficient conditions on a tree $\Upsilon $ for the existence of equivalent 
LUR or strictly convex norms on $\C_0(\Upsilon )$ and for norms with the Kadec 
Property. In particular, we show that for a {\sl finitely branching} tree 
$\Upsilon $ the space $\C_0(\Upsilon )$ admits a Kadec renorming. Since some 
finitely branching trees fail the condition for strictly convex renormability, 
we obtain an example of a Banach space that is Kadec renormable but not 
strictly convexifiable. Consideration of specially tailored examples enables 
us to answer the ``three-space problem'' for strictly convex renorming: there 
exists a Banach space $X$ with a closed subspace $Y$ such that both $Y$ and 
the quotient $X/Y$ admit strictly convex norms, while $X$ does not. We also 
solve a problem about the the property of mid-point locally uniform convexity 
(MLUR), showing that this does not imply LUR renormability. 

In the case of smoothness properties of renormings, we give a necessary and 
sufficient condition for the existence of a Fr\'echet-smooth renorming, which 
turns out (rather surprisingly) to be the same as the condition for LUR 
renormability. We also show that when $\C_0(\Upsilon )$ 
admits a Fr\'echet smooth norm it even admits a $\C^\infty$ norm. By 
considering another specially tailored tree we obtain solution to the quotient 
problem for Fr\'echet differentiable renorming: there is a Banach space $X$ 
with a Fr\'echet differentiable norm and a closed subspace $Y$ such that the 
quotient space $X/Y$ admits no Fr\'echet differentiable renorming. We have few 
results about G\^ateaux differentiability, though we do give an example of a 
space $X$ (as always, of the form $\C_0(\Upsilon )$) which is not strictly 
convexifiable even though it admits an equivalent norm with strictly convex 
dual norm. This seems to be the first example of a space with a G\^ateaux 
smooth norm which is not strictly convexifiable. 

Finally we include some results that are not about renorming theory in the 
strict sense. The theorem of Ekeland and Lebourg already cited in this 
introduction was in fact proved with  hypothesis weaker than the existence of 
a Fr\'echet smooth norm. They needed only a non-trivial Fr\'echet 
differentiable function of bounded support, what is usually referred to as a 
``bump function''. A number of more recent results have been established 
subject to hypotheses of this type, including the interesting 
variational results of Deville, Godefroy and Zizler [\DGZtwo]. Using a 
non-linear version of our technique for the construction of Fr\'echet 
differentiable norms, as presented in [\Haydonfive], we show that 
$\C_0(\Upsilon )$ admits a $\C^\infty$ bump function for every tree $\Upsilon 
$, so that there certainly exist spaces that admit bump functions but no good 
norms. Finally, using another result from [\Haydonfive], we show that 
$\C_0(\Upsilon )$ always admits $\C^\infty$ partitions of unity. 

Much of the early work on this paper was done during the author's stay at 
the \'Equipe d'Analyse, Universit\'e Paris VI in the academic year 
1991--2.  Thanks are due to other members of that team, and especially to 
Gilles Godefroy, for providing a friendly and constantly stimulating working 
environment.
A number of preliminary versions of this work have been circulated since 1991 
and a simple account of the counterexamples to the three-space problem and the 
quotient problem was given in the seminar paper [\Haydonfour]. These examples, 
and the author's early ($\C^1$) bump function result are presented or cited in 
[\DGZ].

\heading 1. Some renorming techniques \endheading

The authors of [\DGZ] record as their Fact II.2.3 the inequality
$$
\phi (x)^2+\phi (y)^2-2\phi (\half(x+y))^2\ge\half(\phi (x)-\phi (y))^2,
$$
valid when $\phi $ is a non-negative convex function. It follows that the 
convergence to 0 of $\phi (x)^2+\phi (x_n)^2-2\phi (\half(x+x_n))^2$ is 
equivalent to the convergence of both $\phi (f_n)$ and $\phi (\half(f_n+f))$ 
to the limit $\phi (f)$. This permits a useful equivalent formulation of the 
definition of local uniform convexity: $\norm{\cdot}$ is LUR if and only if 
the convergence to 0 of $2\norm{x}^2 + 2\norm{x_n}^2-\norm{x+x_n}^2 $ implies 
the norm convergence of $x_n$ to $x$. If, at the end of the last sentence, we 
replace {\sl norm} with {\sl weak} we have the definition of a {\sl weakly} 
locally uniformly convex (wLUR) norm. A different weakening of LUR leads to 
the property known as {\sl midpoint locally uniform convexity} (MLUR): a norm 
is MLUR if the convergence to zero of $\norm{x+h_n}^2 + \norm{x-h_n}^2-
2\norm{x}^2 $ (or, equivalently, the convergence to the limit $\norm{x}$ of 
both $\norm{x+h_n}$ and $\norm{x-h_n}$) implies the norm convergence of $h_n$ 
to $0$. Of course, both wLUR and LUR imply strict convexity, a norm being 
strictly convex if $\norm x=\norm y= \half\norm{x+y}$ implies $x=y$. 

While, as we shall see, it is not strictly speaking a ``convexity'' condition, 
it is convenient to mention here the Kadec property: we say that a norm 
$\norm{\cdot} $ is {\sl Kadec} if the weak and norm topologies coincide on the 
unit sphere $\text{sph}\, X= \{x\in X:\norm x=1\}$.  For Banach spaces $X$ not 
containing $\ell_1$ (and in particular for spaces $X=\C_0(L)$, with $L$ 
scattered) we have the following  sequential characterization: $X$ has the 
Kadec Property if and only if weakly convergence implies norm convergence for 
sequences in the unit sphere. It is not hard to show that a LUR norm is Kadec. 
A striking theorem of Troyanski [\Troy\ or \DGZ, IV.3.6]  states that if a 
Banach space admits both a Kadec renorming and a strictly convex renorming 
then it admits a LUR renorming. 

The Kadec Property is closely related to two other notions, in which the 
linear structure of a Banach space is not involved.  One of these is $\sigma 
$-fragmentability, a property introduced by Jayne, Namioka and Rogers [\JNR].  
We say that the weak topology of a Banach space $X$ is {\sl $\sigma 
$-fragmented} by the norm (or, more briefly, that $X$ is {\sl $\sigma 
$-fragmentable}) if, for every $\epsilon >0$, there is a countable covering 
$(X_n)_{n\in\nats}$ of $X$ such that, for each $n$ and each non-empty subset 
$Y$ of $X_n$, there is a non-empty subset of $Y$ which is relatively open in 
the weak topology and has norm diameter smaller than $\epsilon$. Recent work 
of Namioka and Pol shows that $\sigma$-fragmentability of a Banach space $X$ 
is characterized by a purely topological property of $X$ in the weak topology. 
It is shown in [\JNR] that a Banach space with a Kadec norm is $\sigma 
$-fragmentable, and it is conjectured that any $\sigma $-fragmentable Banach 
space admits a Kadec renorming.  At the end of Section~6 we verify this 
conjecture for spaces $\C_0(\Upsilon )$. We shall say that a locally compact 
space $L$ has the {\sl Namioka Property } $\N\dual$ if, for every Baire space 
$B$ and every function $\psi :B\to \C_0(L)$ which is continuous into the 
topology $\tp$ of pointwise convergence on $L$, the set of points of 
continuity of $\psi $ into the norm topology is dense in $B$. If $L$ is a 
locally compact scattered space and $\C_0(L)$ is $\sigma$-fragmentable then 
$L$ has the Namioka Property. (In general, the Namioka Property is implied by 
$\sigma $-fragmentability of the pointwise topology $\tau_{\text p}$, rather 
than of the weak topology of $\C_0(L)$.)

We now pass from definitions to renorming techniques.
The next lemma contains an idea that is at the heart of many LUR renorming 
proofs. In the following form it may be found in [\DGZ]. 

\proclaim{Lemma 1.1}
Let $X, \norm{\cdot}_0$ be a Banach space, let $I$ be a set and let $(\phi 
_i)_{i\in I}$ and $(\psi _i)_{i\in I}$ be families of non-negative 
convex functions on $X$ which are uniformly bounded on bounded subsets of $X$. 
For $f\in X$, $m\in \nats$ and $i\in I$ define 
$$
\align
\phi (f)&=\sup_{i\in I}\phi _i(f)\\
\theta _{i,m}(f)&=\phi_i(f)^2+2^{-m}\psi _i(f)^2\\
\theta _m(f)&=\sup_{i\in I}\theta _{i,m}(f)\\
\theta (f) &= \norm{f}_0^2+\sum_{m=1}^\infty 2^{-m} \theta _m(f).
\endalign      
$$
If $f$ and $f_n$ in $X$ are such that $\theta (f)+\theta (f_n)-2\theta 
(\half(f+f_n))\to 0$ as $n\to\infty$, then there exists a sequence 
$(i_n)$ in $I$ such that :
\roster
\item $\phi _{i_n}(f)$, $\phi _{i_n}(f_n)$ and $ \phi (f_n)$ all converge to 
$\phi(f)$ as $n\to \infty$;
\item $\psi _{i_n}(f_n)^2+\psi _{i_n}(f)^2-2\psi 
_{i_n}(\half(f_n+f))^2$ tends to 0 as $n\to \infty$.
\endroster
Moreover, there is a norm $\norm\cdot$ on $X$, equivalent to $\norm{\cdot}_0$, 
such that the above conclusion holds whenever $ 2\norm{f}^2+2\norm{f_n}^2-
\norm{f+f_n}^2\to 0. $ \endproclaim

In Section 6 we shall use a variant of this technique, appropriate to the 
construction of Kadec norms, rather than LUR norms. The proof is rather 
easier than that of 1.1.

\proclaim{Proposition 1.2}
Let $X$ be a topological space, let $I$ be a set and let $\phi _i,\psi _i:X\to 
\reals$ ($i\in I$) be uniformly bounded families of lower semicontinuous 
functions.  For each $f\in X$ define 
$$
\align
\phi (f)&=\sup_{i\in I}\phi _i(f)\\
\theta _m(f)&=\sup_{i\in I}\bigl[\phi_i(f)+2^{-m}\psi _i(f)\bigr]\\
\theta (f) &= \sum_{m=1}^\infty 2^{-m} \theta _m(f).
\endalign      
$$
If $(f_n)$ is a sequence which converges to $ f$ in $X$ and is such that 
$\theta (f_n)\to \theta (f)$ then there exists a sequence $(i_n)$ in $I$ such 
that $\phi _{i_n}(f)$ and $\phi _{i_n}(f_n)$ and $\phi (f_n)$ all converge as 
$n\to \infty$ to the limit $\phi (f)$ whilst $\psi _{i_n}(f_n)-\psi _{i_n}(f)$ 
tends to $0$.  If $I$ is equipped with a topology under which it is 
sequentially compact then we may suppose $(i_n)$ to be a convergent sequence 
in this topology.
\endproclaim
\demo{Proof}
For each $m$ the function $\theta _m$ is l.s.c. on $X$.  Since $f_n\to f$ in 
$X$ and $\theta (f_n)\to \theta (f)$ it must therefore be that 
$\theta_m (f_n)\to \theta_m (f)$ for all $m$.  Choose, for each $m$, an 
element $j_m$ of $I$ and a natural number $n_m$ such that
$$
\phi _{j_m}(f)+2^{-m}\psi _{j_m}(f)> \sup_{p\ge n_m}\theta _m(f_p)-2^{-2m}.
$$
By lower semicontinuity we may suppose $n_m$ to be chosen so that 
$\phi _{j_m}(f_p)>\phi _{j_m}(f)-2^{-2m}$ and $\psi _{j_m}(f_p)>\psi 
_{j_m}(f)-2^{-m}$ whenever $p\ge n_m$.  We may also assume that $n_{m+1}>n_m$.

We now let $(m(k))_{k\in\nats}$ be any subsequence of the natural numbers, 
chosen, if we wish, so that $(j_{m(k)})$ converges with respect to a given 
sequentially compact topology on $I$.  Define $i_p=j_{m(k_p)}$ where
$k_p = \max \{k\in \nats:n_{m(k)}\le p\}$. With these definitions, we 
certainly have $p\ge n_{m(k_p)}$ so that
$$
\align
\phi _{i_p}(f)+2^{-m(k_p)}\psi _{i_p}(f) &> \sup_{q\ge n_{m(k_p)}} \theta 
_{m(k_p)}(f_q)-2^{-2m(k_p)}\\
      &\ge \phi _{i_p}(f_p)+2^{-m(k_p)}\psi _{i_p}(f_p)-2^{-2m(k_p)}.
\endalign
$$
Since we also have $\phi _{i_p}(f_p)>\phi _{i_p}(f)-2^{-2m(k_p)}$ and $\psi 
_{i_p}(f_p)>\psi _{i_p}(f)-2^{-m(k_p)}$, it must be that
$$
\align
|\phi_{i_p}(f)-\phi_{i_p}(f_p)| &< 2^{-2m(k_p)+1}\\
|\psi_{i_p}(f)-\psi_{i_p}(f_p)| &< 2^{-m(k_p)+1}.
\endalign
$$
Thus $\phi_{i_p}(f)-\phi_{i_p}(f_p)$ and $\psi_{i_p}(f)-\psi_{i_p}(f_p)$ tend 
to 0 as $p\to \infty$.  We also have the inequalities
$$
\align
\phi _{i_p}(f)+2^{-m(k_p)}\sup_i\inorm{\psi _i} &\ge 
                          \phi _{i_p}(f)+2^{-m(k_p)}\psi _{i_p}(f) \\ 
       &\ge \sup_{q\ge n_{m(k_p)}} \theta _{m(k_p)}(f_q)-2^{-2m(k_p)}\\
       &\ge \limsup_{q\to \infty} \phi (f_q) -2^{-2m(k_p)}\\
       &\ge \liminf_{q\to \infty} \phi (f_q) -2^{-2m(k_p)}\\
       &\ge \phi (f) -2^{-2m(k_p)}\\
       &\ge \phi _{i_p} (f) -2^{-2m(k_p)},
\endalign
$$
which yield the convergence of $\phi _{i_p}(f)$ and $\phi (f_p)$ to the limit 
$\phi (f)$.
\qed \enddemo 

Apart from the above Proposition, the main ingredient in Section~6 will be 
the use of recursion to define norms. There is nothing particularly new in 
this approach, but it may make Section~6 clearer if we give here some more 
elementary examples of recursively defined norms. We start with a construction 
of a LUR renorming of $\C[0,\Omega]$, where $\Omega $ is an ordinal.  The same 
approach may be used in greater generality to give an alternative construction 
of LUR norms on certain spaces having projectional resolutions of the identity 
[\DGZ, Proposition VII.1.6]. It should be emphasized that the existence of 
such norms is well known and that following proposition is included here for 
the (possible) interest of the method, rather than any originality of the 
conclusion. 

\proclaim{Proposition 1.3}
Let $\Omega$ be an ordinal.  There is a unique real-valued function $\Phi$ 
with domain $\C[0,\Omega]\times\{(\alpha ,\gamma ):0\le\alpha \le \gamma\le 
\Omega \}$, satisfying the inequality 
$0\le\Phi (f,\alpha ,\gamma )\le \inorm f$ and the identities
$$
\align
\ \,\Phi (f,\alpha ,\alpha) &= |f(\alpha )| \\
16\Phi (f,\alpha ,\gamma)^2 &=  4\inorm{f\restriction [\alpha ,\gamma ]}^2 
    +f(\alpha )^2 
     +\osc(f\restriction [\alpha ,\gamma ])^2 \\
   +  \sum_{m=1}^\infty 2^{-m}\sup_{\alpha \le \beta <\gamma } 
  &\biggl[(f(\beta +1)-f(\beta ))^2+2^{-m}\Phi (f,\alpha ,\beta )^2
                                     +2^{-m}\Phi (f,\beta +1,\gamma )^2\biggl]
\endalign
$$
for all $f\in \C[0,\Omega )$ and $\alpha <\gamma \le \Omega $. If we define
$\onorm{f}=\Phi (f,0,\Omega )$, then $\onorm{\cdot}$ is a locally uniformly 
convex norm on $\C[0,\Omega ]$ satisfying $\half\inorm{f}\le \onorm{f}\le 
\inorm{f}$.
\endproclaim
\demo{Proof}
As with many results on recursion, it is convenient to give a proof using 
a fixed-point theorem. Let $\Cal D$ be the set of all triples $(f,\alpha 
,\gamma )$ with $f$ in $\C[0,\Omega ]$ and $0\le \alpha \le \gamma \le\Omega$; 
let $\Cal X$ be the space of all functions from $\Cal D$
into the real interval $[0,1]$, which are  positively homogeneous in $f$.
Equip $\Cal X$ with the metric
$$
d(\Phi ,\Psi )=\sup\{|\Phi (f,\alpha ,\gamma )-\Psi (f,\alpha ,\gamma )|:0\le 
\alpha \le \gamma \le\Omega ,\ \inorm f\le 1\}.
$$ 
Define $T:\Cal X\to \Cal X$ by
$$
\align
\ \,(T\Psi) (f,\alpha ,\alpha) &= |f(\alpha )| \\
16(T\Psi) (f,\alpha ,\gamma)^2 &= 4\inorm{f\restriction [\alpha ,\gamma ]}^2 
     +f(\alpha )^2
     +\osc(f\restriction [\alpha ,\gamma ])^2 \\
   +  \sum_{m=1}^\infty 2^{-m}\sup_{\alpha \le \beta <\gamma } 
   &\biggl[(f(\beta +1)-f(\beta ))^2+2^{-m}\Psi (f,\alpha ,\beta )^2
                                      +2^{-m}\Psi (f,\beta +1,\gamma )^2\biggl]
\endalign
$$
The mapping $T$ is a strict contraction on the complete metric space $\Cal X$ 
and so has a unique fixed point $\Phi $. Evidently,  $\onorm{f}=\Phi 
(f,0,\Omega )$ defines a norm on $\C[0,\Omega ]$ satisfying 
$\half\inorm{f}\le \onorm{f}\le \inorm{f}$. We have to show that this norm is 
LUR. Let $f$ be a given element of $\C[0,\Omega ]$ with $\inorm{f}\le 1$ and 
let $\epsilon $ be a positive real number. We shall need an easy result about 
continuous functions on ordinals. 

\proclaim{Fact} {\rm There exists a positive real number $\eta $ such that if 
$0\le \alpha <\gamma\le \Omega $ and $|f(\beta)-f(\beta +1)|<\eta $ for all 
$\alpha \le \beta <\gamma $ then $|f(\alpha )-f(\beta )|<\epsilon $ for all 
$\alpha <\beta \le \gamma $.}\endproclaim

\smallskip
If $|f(\beta )-f(\Omega )|<\half\epsilon $ for all $\beta\in [0,\Omega )$ 
there is nothing to prove. Otherwise, we start a recursion, choosing $\beta _1$ 
to be the greatest ordinal with $|f(\beta_1 )-f(\Omega )|\ge\half\epsilon $. 
Subsequently, if $\Omega >\beta _1>\cdots>\beta _j$ have been defined, then 
either we have $|f(\beta )-f(\beta _j)|<\half\epsilon $ for all $0\le \beta 
<\beta _j$, in which case we stop, of else we take $\beta _{j+1}$ to be the 
greatest ordinal with $0\le \beta_{j+1} <\beta _j$ and $|f(\beta_{j+1} )-
f(\beta _j)|\ge\half\epsilon $.  By well-ordering, the process does stop at 
some stage, after the definition of $\beta _k$, say. We observe that the 
oscillation of $f$ on each of the intervals $[0,\beta _k]$, $(\beta _k,\beta 
_{k-1}] $,\dots $(\beta _1,\Omega ]$ is smaller than $\epsilon $ and that 
$f(\beta _j)\ne f(\beta _j+1)$ for each $j$. We set $\eta =\min\{|f(\beta _j)-
f(\beta _j+1)|:1\le j\le k\}$. If $\alpha <\gamma $ are such that the interval 
$[\alpha ,\gamma )$ contains no $\beta $ with $|f(\beta )-f(\beta +1)|\ge \eta 
$, then $[\alpha ,\gamma ]$ must be contained in one of the intervals 
$[0,\beta _k]$, $(\beta _k,\beta _{k-1}] $,\dots $(\beta _1,\Omega ]$.  Thus 
the oscillation of $f$ on $[\alpha ,\gamma ]$ is smaller than $\epsilon $, as 
required. 

\smallskip
We now define $m(\alpha ,\gamma )$ to be the (finite!) number of $\beta $ such 
that $\alpha \le \beta <\gamma $ and $|f(\beta )-f(\beta +1)|\ge \eta $. We 
shall use induction on $m(\alpha ,\gamma )$ to establish the following 
assertion, which will obviously complete the proof of 1.3.

\proclaim{Claim}{\rm If $f_n\in \C[0,\Omega ]$ are such that $2\Phi (f_n,\alpha 
,\gamma )^2+2\Phi (f,\alpha ,\gamma )^2-\Phi (f+f_n,\alpha ,\gamma )^2\to 0$ 
as $n\to \infty$ then $\limsup\inorm{(f-f_n)\restriction [\alpha ,\gamma 
]}<4\epsilon $.  }\endproclaim

\smallskip
If $m(\alpha ,\gamma )=0$ then, by the way in which we chose $\eta $, the 
oscillation of $f$ on $[\alpha ,\gamma ]$ is smaller than $\epsilon $, so 
that $\osc(f\restriction [\alpha ,\gamma ]) <\epsilon  $.
>From our hypothesis about $f_n$ and convexity arguments, we see that 
$f_n(\alpha)$ tends to $f(\alpha)$ and $\osc(f_n\restriction[\alpha,\gamma])$ 
tends to $\osc(f\restriction [\alpha ,\gamma ]) $ as $n\to \infty$.  Thus, for 
all large enough $n$ we have $|f(\alpha )-f_n(\alpha )|<\epsilon $ and 
$\osc(f_n\restriction[\alpha,\gamma])<2\epsilon $, giving $\inorm{(f-
f_n)\restriction [\alpha ,\gamma ]}<4\epsilon $.

We now assume inductively that if $\beta $ and $\delta $ are such that 
$m(\beta ,\delta )<m(\alpha ,\gamma )$ and $(g_n)$ is a sequence such that 
$2\Phi (g_n,\beta ,\delta )^2+2\Phi 
(f,\beta ,\delta )^2-\Phi (f+g_n,\beta ,\delta )^2\to 0$, then 
$\limsup\inorm{(f-g_n)\restriction [\beta ,\delta ]}<4\epsilon $. If the 
assertion we are trying to prove is false, then we may, by passing to a 
subsequence, assume that $\inorm{(f-f_n)\restriction[\alpha ,\gamma ]}\ge 
4\epsilon $ for all $n$.  We are finally ready to apply Proposition~1.1 with
$$
\align
\phi_\beta (f) &= |f(\beta +1)-f(\beta )|,\\
\psi_\beta (f)  &= \left(\Phi (f,\alpha ,\beta )^2+\Phi (f,\beta +1,\gamma 
)^2\right)^{\frac12}.
\endalign
$$
There exists a sequence $(\beta _n)$ such that the conclusions of 1.1 hold, in 
particular,
$$
\align
|f(\beta_n +1)-f(\beta_n )|&\to \sup_{\alpha \le \beta <\gamma }
                                                      |f(\beta +1)-f(\beta)|\\
2\Phi (f_n,\alpha ,\beta _n )^2+2\Phi 
(f,\alpha ,\beta _n )^2-\Phi (f+f_n,\alpha ,\beta _n )^2&\to 0    \\
2\Phi (f_n,\beta_n ,\delta )^2+2\Phi 
(f,\beta_n ,\delta )^2-\Phi (f+f_n,\beta_n ,\delta )^2&\to 0        .
\endalign
$$
For all large enough $n$, it must be that $\beta _n$ is in the finite set of 
$\beta $ with $|f(\beta )-f(\beta +1)|\ge \eta $ and, for a suitable 
subsequence,  
$\beta _{n_k}$ will take the same value $\beta $. Since we have $m(\alpha 
,\beta )<m(\alpha ,\gamma )$ and $m(\beta ,\gamma )<m(\alpha ,\gamma )$  our 
inductive hypothesis is applicable, giving 
$$
\align\limsup\inorm{(f-f_{n_k}\restriction[\alpha ,\beta ]}&<4\epsilon \\
\limsup\inorm{(f-f_{n_k}\restriction[\beta, \gamma ]}&<4\epsilon ,\\
\intertext{whence}
\limsup\inorm{(f-f_{n_k}\restriction[\alpha, \gamma ]}&<4\epsilon ,
\endalign
$$
which completes the proof.\qed
\enddemo

It may amuse the reader to note that the same approach could be used to 
establish the local uniform convexity of Day's norm on $c_0(\Gamma )$.  We 
recall that we define $\norm{\cdot}_{\text{Day}}$ by
$$
\norm{f}_{\text{Day}}^2 = \sup \sum_{n=1}^\infty 2^{-n}f(\gamma_n)^2,
$$
where the supremum is taken over all sequences $(\gamma _n)$ of distinct 
elements of $\Gamma $. One could also regard $\norm{f}_{\text{Day}}$ as being 
$\Phi (f;\Gamma )$, where $\Phi $ is the unique function defined on the 
Cartesian product $c_0(\Gamma )\times \Cal P(\Gamma )$ which satisfies the 
inequality $0\le \Phi (f;\Delta )\le \inorm f$ as well as the identities 
$$
\align
\Phi (f;\emptyset) &= 0\\
\Phi (f;\Delta )^2 &=  \sum_{m=1}^\infty2^{-m}\sup_{t\in \Delta }
\biggl[\half f(t)^2+ ({\textstyle\frac23})^m\Phi (f;\Delta \less\{t\})^2\biggr]
\endalign
$$
for $f\in c_0(\Gamma )$ and $\emptyset \ne \Delta \subseteq \Gamma $. 
Simplifying the proof of 1.3, we see that an argument using 1.1 and proceeding 
by induction on the number of $\gamma$ for which $|f(\gamma )|\ge \epsilon $ 
enables to show that $\norm\cdot_{\text{Day}}$ is LUR.

\smallskip
We now pass to questions of smoothness.
A real-valued function $\phi $, defined on an open subset $U$ of a Banach 
space $X$ is said to be {\sl G\^ateaux-differentiable} at $x\in U$ if there 
exists an element $\eta $ of the dual space $X\dual$ such that, for all $h\in 
X$, $t^{-1}(\phi (x+th)-\phi (x))\to \<\eta ,h\>$ as $t\to 0$. It is usual to 
write $D\phi (x)$ for this element $\eta$. In the case where $\psi $ is a 
norm, there is a simple criterion: $\norm{\cdot}$ is  G\^ateaux-differentiable 
at $x$ if and only if there is a {\sl unique} $\xi \in X\dual$ satisfying 
$\norm \xi \dual =1$, $\<\xi ,x\>=\norm{x}$; moreover, $D\norm{\cdot}(x)$ is 
this element $\xi $. For economy of notation, we shall write $x\dual$ for 
$D\norm{\cdot}(x)$. A norm is said to be {\sl Gateaux-smooth} (or simply {\sl 
smooth}) if it Gateaux-differentiable at all points except 0.

By definition, a function $\phi :U\to \reals$ is {\sl 
Fr\'echet-differentiable} at $x$ if it is G\^ateaux-differentiable and the 
convergence of $t^{-1}(\phi (x+th)-\phi (x))$ to $\<D\phi (x),h\>$ as $t\to 0$ 
is uniform over $h$ in the unit ball of $X$. In the case of a norm, we have 
the important criterion of \v Smulyan [\Smulyan\ or \DGZ, I.1.4]: a norm 
$\norm{\cdot}$ is Fr\'echet-differentiable at $x$ if and only if it is 
G\^ateaux-differentiable at $x$ and we have $\norm{\xi _n-x\dual}\to 0$, 
whenever $(\xi _n)$ is a sequence in $\ball X\dual$ with $\<\xi _n,x\>\to 
\norm{x}$.

For the construction of Fr\'echet-smooth norms our main method will be one
based on Talagrand's construction [\Tal] for the space $\C[0,\Omega ]$. When 
modified as in [\Haydontwo] and [\Haydonfive], this approach even yields 
norms which are infinitely differentiable (except at 0 of course). Given a 
locally compact space $L$ and a set $M$, we shall say that a bounded linear 
operator $\C_0(L)\to c_0(L\times M)$ is a {\sl Talagrand operator} for $L$ if 
for every non-zero $f$ in $\C_0(L)$, there exist $t\in L$ and $u\in M$ with 
$|f(t)|=\inorm f$ and $(Tf)(t,u)\ne 0$. 

\proclaim{Proposition 1.4 {\smc(Theorem 1 of [\Haydonfive])}}
Let $L$ be a set and let $U(L)$ be the subset of the direct sum $\ell_\infty(L)
\oplus c_0(L)$ consisting of all pairs $(f,x)$ such that $\inorm{f}$ and 
$\inorm{x}$ are both strictly less than $\inorm{|f|+\half |x|}$. The space 
$\ell_\infty(L) \oplus c_0(L)$ admits an equivalent norm $\norm{\cdot}$ with 
the following properties: 
\roster
\item $\norm{\cdot}$ is a lattice norm, in the sense that $\norm{(g,y)}\le 
\norm{(f,x)}$ whenever $|g|\le |f|$ and $|y|\le |x|$;
\item $\norm{\cdot}$ is infinitely differentiable on the open set $U(L)$;
\item locally on $U(L)$, $\norm{(f,x)}$ depends on only finitely many non-zero 
coordinates; that is to say, for each $(f^0,x^0)\in U(L)$ there is 
a finite $N\subseteq L$ and an open neighbourhood $V$ of $(f^0,x^0)$ in $U(L)$, 
such that for $(f,x)\in V$ the norm $\norm{(f,x)}$ is determined by the values 
of $f_t$ and $x_t$ with $t\in N$ and such that $f_t\ne0$, $x_t\ne 0$ for all 
such $(f,x)$ and $t$. 
\endroster \endproclaim 

The following Corollary appears in both [\Haydontwo] and [\Haydonfive].

\proclaim{Corollary 1.5 }
Let $L$ be a locally compact space which admits a Talagrand operator.  Then
$\C_0(L)$ admits a $\C^\infty$ renorming.
\endproclaim
\demo{Proof}
Let $T:\C_0(L)\to c_0(L\times M)$ be a Talagrand operator, normalized so that 
$\inorm{Tf}\le \half \inorm f$ for all $f$ and let $S:\C_0(L)\to 
\ell_\infty(L\times M)$ be defined by $(Sf)(t,u)=f(t)$.  The pair $(Sf,Tf)$ is 
in $U(L\times M)$ whenever $f$ is a non-zero element of $\C_0(L)$ and so  the 
norm on $\C_0(L)$ defined by $\norm {f} = \norm{(Sf,Tf)}$ is infinitely 
differentiable except at 0. \qed \enddemo 

\heading
2. Some preliminaries about trees
\endheading

By definition, a tree is a partially ordered set $(\Upsilon ,\trle)$ with the 
property that for every $t\in \Upsilon $ the set $\{s\in\Upsilon :s\trle t\}$ 
is well-ordered by $\trle$. In any tree, we use normal interval notation, so 
that, for instance, $(s,u]=\{t\in\Upsilon :s\trl t\trle u\}$. For convenience 
of notation, we introduce two \lq\lq imaginary'' elements, not in $\Upsilon $,
denoted $0$ and $\infty$, and having the property that $0\prec t\prec\infty$ 
for all $t\in \Upsilon $. This allows us to extend our interval notation to 
include expressions like $(0,t]$ and $[t,\infty)$. Note that, by definition, 
each $(0,t]$ is well-ordered, but that $[t,\infty)$ need not be. For each 
$t\in \Upsilon $ there is a unique ordinal $r(t)$ with the same order type as 
$(0,t)$. Fremlin [\Fremlin] says that a tree is {\sl Hausdorff} if, whenever 
$r(t)$ is a limit ordinal and $(0,t')=(0,t)$, we necessarily have $t=t'$. Such 
a tree may be equipped with a locally compact, and Hausdorff (!), topology 
which may be characterized as the coarsest for which all intervals $(0,t]$ are 
open and closed. We shall consider only trees that are Hausdorff. We shall 
be studying norms on the space $\C_0(\Upsilon )$ of real-valued functions $f$ 
on $\Upsilon $, which are continuous for the locally compact topology and are 
such that, for all $\epsilon >0$ the set $\{t\in \Upsilon :|f(t)|\ge \epsilon 
\}$ is compact for that topology. It will occasionally be convenient to give 
meaning to the expressions $f(0)$ and $f(\infty)$ when $f\in \C_0(\Upsilon )$, 
taking both of them to equal $0$. It is useful to note that the space 
$\C_0(\Upsilon )$ is the closed linear span in $\ell_\infty(\Upsilon )$ of the 
indicator functions $\indic_{(0,t]}$ $(t\in \Upsilon )$. Many results about 
general elements of $\C_0(\Upsilon )$ may be established easily using uniform
approximation by linear combinations of these indicators. Here are two 
examples. 

\proclaim{Lemma 2.1}
Let $\Upsilon $ be a tree, let $f$ be in $\C_0(\Upsilon )$ and let $\delta$ 
be a positive real number. There exists $k\in\nats$ such that if $H$ is a 
subset of $\Upsilon $ whose elements are pairwise incomparable then $|f(t) 
|\ge \delta $ for at most $k$ elements $t$ of $H$. 
\endproclaim

\proclaim{Lemma 2.2}
Let $\Upsilon $ be a tree, let $f$ be in $\C_0(\Upsilon )$ and let $\delta$ 
be a positive real number. For all but finitely many $s\in \Upsilon $, there 
exists $t\in s^+$ such that $|f(t)-f(s)|<\delta $, while $|f(u)|<\delta $ for 
all $v\in [s,\infty)\less[t,\infty)$.
\endproclaim
                                                  
For any tree $\Upsilon $ and any $t\in \Upsilon $ we write $t^+$ for the set 
of all {\sl immediate successors } of $t$ in $\Upsilon $, that is to say that 
$t^+$ contains those elements $u$ with the property that $s\prec u$ iff 
$s\trle t$. By convention, we write $0^+$ for the set of minimal elements of 
$\Upsilon $ and define $\Upsilon ^+= \bigcup_{t\in\Upsilon \cup\{0\}} t^+$. 
The elements of $\Upsilon ^+$ are exactly those $t$ for which $r(t)$ is 0 or a 
limit ordinal; they may also be characterized as the isolated points of 
$\Upsilon $ for its locally compact topology. For $u\in \Upsilon ^+$ we  write 
$u^-$ for the unique $t\in \Upsilon\cup\{0\} $ such that $u\in t^+$. We say 
that $\Upsilon $ is {\sl (in-) finitely branching } if $t^+$ is (in-) finite 
for all $t\in \Upsilon $; we say that $\Upsilon $ is a dyadic tree if each 
$t^+$ contains two elements. A subset $\Gamma $ of a tree $\Upsilon $ is said 
to be {\sl ever-branching} if for every $t\in \Gamma $ the intersection 
$\Gamma \cap [t,\infty)$ is not totally ordered. 

The {\sl height } of a tree $\Upsilon $ is defined to be $ht(\Upsilon 
)=\sup\{r(t)+1:t\in \Upsilon \}$; a {\sl branch of } of $\Upsilon $ is a 
maximal totally ordered subset; we say that $\Upsilon $ is a {\sl full 
tree} of height $\alpha $ if every branch has order-type $\alpha $. The tree 
considered in [\Haydon] was a full uncountably branching tree of height 
$\omega _1$. 

We say that two elements $s,t$ of a tree $\Upsilon $ are {\sl incomparable} if 
neither $s\trle t$ nor $t\trle s$ holds. A subset of $\Upsilon $ whose 
elements are pairwise incomparable  is called an {\sl antichain}. For a subset 
$S$ of $\Upsilon $ we define $\min S$ to be the set of elements of $S$ that 
are minimal for the tree order $\trle$. If $s\in S$ then $(0,s]$ contains a 
member of $\min S$ (because $(0,s]$ is well-ordered). We write $\max S$ for 
the set of maximal elements of a subset $S$ of $\Upsilon$. In general, 
non-empty subsets of $\Upsilon $ need not have maximal elements, but if $S$ is 
a {\sl compact} subset of $\Upsilon $ then $\max S$ is finite and $S\subseteq 
\bigcup_{s\in\max S}(0,s]$ (since the sets $(0,s]$ $(s\in S)$ form an open 
cover of $S$). Of course, for any $S\subseteq\Upsilon $, the sets $\min S$ and 
$\max S$ are antichains. 

There are other interesting topologies with which a tree may be equipped, as 
well as the locally compact topology we have been looking at so far. Logicians 
are often concerned with what is sometimes called the {\sl``forcing 
topology''}, whose basic open sets are of the form $[t,\infty)$. This 
satisfies the $\text T_0$ separation axiom, but (except in the trivial case 
where $\Upsilon $ is an antichain) not the $\text T_1$ axiom. We shall be 
interested in a variant of this topology, which will here be called the {\sl 
reverse topology}, and which we define to be the coarsest topology for which 
all the subsets $[t,\infty)$ are open and closed. A base of neighbourhoods of 
$t$ for the reverse topology consists of the sets $$ 
[t,\infty)\less\bigcup_{u\in F} [u,\infty), $$ with $F$ a finite subset of 
$t^+$. Evidently $t$ is an isolated point for this topology if and only if 
$t^+$ is finite, so that $\Upsilon $ is reverse-discrete if and only if 
$\Upsilon $ is  finitely branching. The following observation goes some way 
towards explaining why the reverse topology is of interest. 

\proclaim{Proposition 2.3}
For any tree $\Upsilon $, the map $u\mapsto \indic_{(0,u]}$ is a homeomorphic 
embedding from $\Upsilon $ in the reverse topology into $\C_0(\Upsilon )$ in 
the topology $\tp$  of pointwise convergence. 
\endproclaim
\demo{Proof} Restricted to the set of indicator functions, the topology $\tp$  
may be characterized as the coarsest for which all sets $\{f: f(t)=1\}$ 
$(t\in\Upsilon $) are open and closed. Since $\indic_u(t)=1$ if and only if 
$u\in[t,\infty)$, we see that this corresponds exactly to our definition of 
the reverse topology. \qed\enddemo

\proclaim{Corollary 2.4}
If $\Upsilon $ is an infinitely branching tree which is a Baire space for the 
reverse topology then $\Upsilon $ does not have the Namioka property.
\endproclaim
\demo{Proof}
We have already noted that an infinitely branching tree has no 
reverse-isolated points, so that $u\mapsto \indic_{(0,u]}$ has no points of 
continuity from the reverse topology into the norm topology of $\C_0(\Upsilon 
)$. On the other hand, this map is continuous into $\tp$ and, by hypothesis, 
$\Upsilon $ is reverse-Baire.
\qed\enddemo

The easy argument at the start of [\Haydon] shows that a full tree of height 
$\omega _1$ is reverse Baire. However, examples exist [\Tod\ or \Haydonfour] 
of infinitely-branching reverse-Baire trees in which every branch is 
countable. It is shown in [\Haydonfour] that infinitely branching Baire trees 
have all the bad properties of the tree considered in [\Haydon]. Of course, a 
consequence of this is that if we want to study the fine structure of trees in 
renorming theory we have to look at trees which are not Baire for the reverse 
topology. 

As already mentioned in the Introduction, we shall obtain necessary and 
sufficient conditions for the existence on $\C_0(\Upsilon )$ of various good 
renormings. These conditions will all be expressed in terms of increasing 
real-valued functions on $\Upsilon $. A function $\rho :\Upsilon \to\reals$ is 
said to be {\sl increasing} if $s\trle t\implies \rho(s)\le \rho(t)$ and to be 
{\sl strictly} increasing if $s\trl t\implies \rho(s)< \rho(t)$. A tree is 
said to be $\reals$-{\sl embeddable} if it admits a strictly increasing 
real-valued function. A tree which admits a strictly increasing rational-valued 
function may be said to be $\rats$-embeddable, but it is more usual to call 
such a tree {\sl special}. An equivalent definition is to say that a tree is 
special if it is a countable union of antichains. If Baire trees are very bad 
from the point of view of renorming theory, then special trees are very good. 
This may be readily deduced from known results. 

\proclaim{Proposition 2.5} Let $L$ be a locally compact space which is a 
countable union of closed subsets $L_n$ $(n\in \nats)$. If, for each $n$, 
$\C_0(L_n)$ admits an equivalent norm which is LUR and has LUR dual norm, then 
the same is true for $\C_0(L)$.
\endproclaim
\demo{Proof}
For each $n$, let $\norm{\cdot}_n$ be an equivalent norm on $\C_0(L_n)$, which 
is LUR, with LUR dual norm.  We may suppose that $\norm{\cdot}_n\le 
\norm{\cdot}_\infty$.  The $\ell_2$ direct sum $Y=\left(\bigoplus 
(\C_0(L_n),\norm{\cdot}_n)\right)_2$, may be equipped with its natural
norm $\norm{\cdot}_Y$, where
$$
\norm{(g_n)}_Y^2 = \sum_n \norm{g_n}_n^2.
$$
This norm is LUR and has LUR dual norm.  We define a bounded linear operator 
$T:\C_0(L)\to Y$ by $(Tf)_n = f\restriction_{L_n}$ and note that the image of 
the transpose operator $T\dual$ contains each Dirac functional $\delta _t$ 
$(t\in L)$.  Thus $T^{\textstyle **}$ is injective and Theorem VII.2.6 of 
[\DGZ] may be applied to give the required renorming of $\C_0(L)$.\qed
\enddemo

\proclaim{Corollary 2.6} If $\Upsilon $ is a special tree then $\C_0(\Upsilon 
)$ admits an equivalent norm which is LUR and has LUR dual norm.
\endproclaim
\demo{Proof}
If $A$ is an antichain in $\Upsilon $ then $A$ is closed and discrete in the 
locally compact topology of $\Upsilon $. Because $A$ is discrete,
$\C_0(A)=c_0(A)$, a space with a renorming which is LUR and LUR*. Thus,
if $\Upsilon $ is expressible as a countable union of antichains, Proposition 
2.5 is applicable.
\qed\enddemo

Although less nice than special trees, $\reals$-embeddable trees still have 
some good properties.

\proclaim{Proposition 2.7}
If $\Upsilon $ is $\reals$-embeddable, there exists a bounded linear injection 
from $\C_0(\Upsilon )$ into $c_0(\Upsilon )$.  Hence  $\C_0(\Upsilon )$ admits 
an equivalent strictly convex norm.
\endproclaim
\demo{Proof}
Let $\rho :\Upsilon \to \reals$ be strictly increasing.  By replacing $\rho $,
if necessary, with $\text e^\rho /(1+\text e^\rho )$, we may supose that $\rho 
$ takes values in the real interval $(0,1)$. We define 
$$
(Rf)(t) = \left\{\matrix 
                   \format \c & \l \\  
(\rho (t)-\rho (t^-))f(t) & \qquad\qquad \text{if }t\in \Upsilon ^+\\
                   0  & \qquad\qquad \text{otherwise},
          \endmatrix \right.
$$                                                 
and note that $R$ is a bounded linear operator from $\C_0(\Upsilon )$ into 
$\ell_\infty(\Upsilon )$. If $f$ is in the kernel of $R$, then $f(t)=0$ for 
all $t\in \Upsilon ^+$ because $\rho (t)\ne \rho (t^-)$ for all such $t$.  
Since $\Upsilon ^+$ is dense in $\Upsilon $, we deduce that $f$ is everywhere, 
and hence that $R$ is injective.
 
To show that $R$ takes values in $c_0(\Upsilon )$ it is enough to show that 
$R\indic_{(0,u]}$ is in $c_0(\Upsilon )$ for all $u\in \Upsilon $. This may be 
done by a calculation since 
$$
\align
\sum_{t\in \Upsilon }| (Rf)(t)| &= \sum_{t\in (0,u]\cap \Upsilon ^+} (\rho 
(t)-\rho (t^-))\\
       &\le \rho (u)\le 1,
\endalign
$$
since we are adding up the jumps of the non-negative increasing function $\rho 
$ along the totally ordered set $(0,u]$. What we have shown is that 
$R\indic_{(0,u]}$ is even in $\ell_1(\Upsilon )$.
\qed\enddemo

Special trees and $\reals$-embeddable trees are well-established in the 
literature of Logic and Set Theory. In the next section we shall pass on to a 
definition that is motivated by our applications in renorming.

\heading 3. Good points, bad points and $\mu $-functions \endheading

Given a tree $\Upsilon $ and an increasing function $\rho :\Upsilon \to 
\reals$, we shall say that an element $t$ of $\Upsilon $ is a {\sl good point} 
for $\rho $ if there is a finite subset $F$ of $t^+$ such that $\inf_{u\in 
t^+\less F}\rho (u)>\rho (t)$. Provided there is no ambiguity about which is 
the function $\rho $ under consideration, we shall write $F_t=\{u\in t^+:\rho 
(u)=\rho (t)\}$ and $\delta _t=\inf\{\rho (u)-\rho (t):u\in t^+\less F_t\}.$ 
An element of $\Upsilon $ which is not a good point for $\rho $ will be called 
a {\sl bad point} for $\rho $. 

We can get an idea of the way in which good points and bad points are going to 
enter into the theory by establishing straightaway some necessary conditions 
for the existence on $\C_0(\Upsilon )$ of an equivalent norm with the Kadec 
property and for the existence of an equivalent strictly convex norm. 
Whenever we have a renorming $\norm{\cdot}$ of $\C_0(\Upsilon )$, we 
may define an increasing function $\mu :\Upsilon \to\reals$ by
$$
\mu (t)=\inf\{\norm{f}: f=\indic_{(0,t]}+f':\supp{f'}\subseteq(t,\infty)\}.
$$
The function $\mu $ defined in the lemma that follows, and variants of it,  
will be used throughout the present work.

\proclaim{Lemma 3.1}
Let $\norm{\cdot}$ be an equivalent norm on $\C_0(\Upsilon) $ and let $\mu $ 
be the increasing function defined by
$$ 
\mu (t)=\inf\{\norm{f}: f=\indic_{(0,t]}+f'\text{ with }
\supp{f'}\subseteq(t,\infty)\}. 
$$ 
If $t$ is a bad point for $\mu $ then $\norm{\indic_{(0,t]}}=\mu (t)$. 
\endproclaim 
\demo{Proof}
It is clear from the definition of $\mu $ that $\norm{\indic_{(0,t]}}
\ge\mu (t)$.
Badness implies that there is a sequence $(u_n)$ of distinct elements of $t^+$ 
such that $\mu (u_n)\to \mu (t)$ as $n\to\infty$. Thus there are functions 
$f_n$ of the form $f_n= \indic_{(0,u_n]}+f'_n$, with $\supp f'_n\subseteq 
(u_n,\infty)$, such that $\norm{f_n}\to \mu (t)$. Now $(f_n)$ is a 
norm-bounded sequence and converges pointwise on $\Upsilon $ to 
$\indic_{(0,t]}$. Thus this sequence converges weakly to $\indic_{(0,t]}$, 
which shows that $\norm{\indic_{(0,t]}}\le\lim\norm{f_n}=\mu (t)$.
\qed\enddemo

\proclaim{Proposition 3.2}
If $\C_0(\Upsilon )$ admits an equivalent Kadec norm then the function $\mu $ 
associated with such a norm has no bad points.
\endproclaim
\demo{Proof}
We observe that if $t$ is a bad 
point then the sequence $(f_n)$ defined in the proof of Lemma 3.1 converges 
weakly to $\indic_{(0,t]}$ and satisfies $\norm{f_n}\to \norm{\indic_{(0,t]}}$,
but does not converge in norm, since $\inorm{f_n-\indic_{(0,t]}}\ge 1$ for all 
$n$. Since this contradicts the Kadec property of the norm, $\mu $ can in 
fact have no bad points. 
\qed\enddemo

\proclaim{Proposition 3.3}
Let $\Upsilon $ be a tree and suppose that $\C_0(\Upsilon )$ admits a strictly 
convex renorming. Then the associated function $\mu $ has the following 
property: for every $s \in \Upsilon $ there is at most one bad point $t$ 
satisfying $t\trge s$ and $\mu (t)=\mu (s)$. In particular, $\mu $ is
strictly increasing on the set of its bad points.
\endproclaim
\demo{Proof} 
Suppose that $s$ is in $\Upsilon $ and that $t,u$ are bad points with 
$s\trle t$, $s\trle u$ and $\mu (s)=\mu (t)=\mu (u)=\alpha $. Then we have
$$
\norm{\indic_{(0,t]}}=\alpha , \quad\norm{\indic_{(0,u]}}=\alpha 
$$
by Lemma 3.1, while
$$
\norm{\half(\indic_{(0,t]}+\indic_{(0,u]})}=
\norm{\indic_{(0,s]}+\half(\indic_{(s,t]}+\indic_{(s,u]})}\ge \alpha ,
$$
by the definition of $\mu(s) $. The strict convexity of the norm implies that 
we must have $\indic_{(0,t]}=\indic_{(0,u]}$ so that $t=u$. \qed \enddemo

The next proposition gives a second property possessed by the $\mu $-function 
associated with a strictly convex norm. If $\rho :\Upsilon \to \reals$ is 
increasing, we shall say that a point $u$ of $\Upsilon $ is a {\sl fan point} 
for $\rho $ if $u$ is a member of an ever-branching subset on which $\rho $ is 
constant. There is a useful construction that one may carry out in this case.
Suppose that $\rho $ takes the constant value $\alpha $ on some ever-branching 
subset $T$ of $\Upsilon $. Since $T $ is ever-branching, we can choose for any 
$u$ in $T $, two incomparable elements of $(u,\infty)\cap T$, $u_0$ and $u_1$ 
say. Similarly we may choose incomparable elements $u_{00}$ and $u_{01}$ of 
$(u_0,\infty)\cap T$, and it is possible to continue in such a way as to embed 
a full dyadic tree of height $\omega$ in $[u,\infty)\cap T$, the nodes being 
labelled $u_\sigma $ ($\sigma \in \{0,1\}^n$, $n\in \nats$). We define a 
function $\phi_u\in \C_0(\Upsilon )$ by 
$$ 
\phi_u =  \half(\indic_{(u,u_0]}+\indic_{(u,u_1]})+ \sum\Sb n\ge 1\\
\sigma \in \{0,1\}^n 
\endSb 2^{-n-1} (\indic_{(u_\sigma ,u_{\sigma 0}]}+\indic_{(u_\sigma 
,u_{\sigma 1}]}). 
$$ 
Notice that $\phi _u$ takes the value $\half$ at $u_0$ 
and $u_1$, the value $\quarter$ at $u_{00},u_{01},u_{10},u_{11}$, and so on. 
We shall call $\phi _u$ a {\sl fan function} for $\rho $ at $u$.

\proclaim{Proposition 3.4}
If $\Upsilon $ is a tree and $\norm{\cdot}$ is a strictly convex renorming of 
$\C_0(\Upsilon )$ then the corresponding function $\mu :\Upsilon \to\reals$ is 
constant on no ever-branching subset of $\Upsilon $.
\endproclaim
\demo{Proof}
Suppose if possible that $\mu $ takes the constant value $\alpha $ on some 
ever-branching subset $T$. We fix $u$ in $T$ and introduce $u_\sigma $ 
$(\sigma \in\bigcup_{n\in\nats}\{0,1\}^n)$ as above.  We may choose, for 
each $n$ and each $\sigma \in \{0,1\}^n$, an element $\psi _\sigma $ of $Y$ 
with $\supp \psi _\sigma \subseteq (u_\sigma ,\infty)$ and 
$\norm{\indic_{(0,u_\sigma ]}+\psi _\sigma }\le\alpha+2^{-n}$. Notice that 
$\inorm{\psi _\sigma }\le \inorm{\indic_{(u,u_\sigma ]}+\psi _\sigma }$ 
because $\psi _\sigma $ is supported on $(u_\sigma ,\infty)$ and that 
$\inorm{\psi _\sigma }$ is thus at most $M(\alpha +2^{-n})$. 

For each $n$, we set  
$$
y_n = 2^{-n}\sum_{\sigma \in\{0,1\}^n} \biggl[\indic_{(0,u_\sigma ]}+\psi 
_\sigma 
            \biggr]
$$
and note that $\norm{y_n}\le \alpha +2^{-n}$. Moreover, if $\phi _u$ is the 
fan function introduced above, 
$\indic_{(0,u]}+\phi _u-y_n$ is supported on the union of the disjoint sets 
$(u_\sigma ,\infty)$  ($\sigma \in \{0,1\}^n$) and we have
$$
\align
\inorm{\indic_{(0,u]}+\phi _u-y_n}&\le \max_{\sigma \in\{0,1\}^n}
\inorm{\phi _u\restriction_{(u_\sigma ,\infty)}} + 
2^{-n}\inorm{\psi_\sigma }\\
&\le 2^{-n-1}+2^{-n}M(\alpha +2^{-n}).
\endalign
$$
Thus $\indic_{(0,u]}+\phi _u$ is the norm limit of $y_n$ and so satisfies
$\norm{\indic_{(0,u]}+\phi _u}=\alpha $.  

Now the same argument may be applied with $u$ replaced in turn by $u_0$ and by 
$u_1$.  We find that 
$$
\norm{\indic_{(0,u_0]}+\phi _{u_0}}=\norm{\indic_{(0,u_1]}+\phi_{u_1}}=\alpha,
$$
where $\phi _{u_i}$ are defined in the obvious way.  This contradicts strict 
convexity of the norm $\norm{\cdot}$ since 
$$
\indic_{(0,u]}+\phi _u = \half\biggl(\indic_{(0,u_0]}+\phi _{u_0}+
 \indic_{(0,u_1]}+\phi _{u_1} \biggr).
$$
\qed\enddemo

The propositions that we have just proved give a good idea of why bad points 
and fan points are relevant to convexity properties of equivalent norms on 
$\C_0(\Upsilon )$.  The connection with smoothness is suggested by the next 
one. It is convenient to introduce here a more general notion of $\mu 
$-function.  If $t$ is in $\Upsilon $ we shall, as in [\Haydonthree] write 
$C_t$ for the set of $f\in\C_0(\Upsilon )$ such that $\supp f=(0,t]$ and such 
that $f\restriction (0,t]$ is increasing. For $f\in C_t$ and $u\in [t,\infty)$ 
we set 
$$
\mu(f,u)=\inf\{\norm{f+f(t)\indic_{(t,u]}+g}: g\in \C_0(\Upsilon ), \ \supp 
g\subseteq (u,\infty)\}.
$$
Notice that, if $f\in C_t$ and $t\trle u \trle v$ we can rewrite $\mu(f,v)$ as
$\mu(g,v)$ with $g=  f+f(t)\indic_{(t,u]}\in C_u$. 

\proclaim{Proposition 3.5}
Let $\Upsilon $ be a tree, let $t$ be in $\Upsilon $ and let $f$ be in $C_t$.
Suppose that there exists a function $\hat f$ which attains the infimum in the 
definition of $\mu (f,t)$, that is to say that $\supp (\hat f-f)\subseteq 
(t,\infty)$ and $\norm {\hat f} =\mu (f,t)$. If the norm $\norm\cdot$ is 
Gateaux-differentiable at $\hat f$ (with differential $\hat f\dual$), then 
necessarily $ \<\hat f,h\>=0$ for all $h$ with $\supp h\subseteq (t,\infty)$.
This situation occurs in particular for $\hat f=f$ if $t$ is a bad point for 
$\mu (f,\cdot)$ and for $\hat f=f+f(t)\phi _t$ if $t$ is a fan point for this 
function.
\endproclaim
\demo{Proof}
The assumption about $\hat f$ is that, for all $h$ with $\supp h\subseteq 
(t,\infty)$, the function $\lambda \mapsto \norm{\hat f+\lambda h}:\reals\to 
\reals$ has a minimum at $\lambda =0$.  By Gateaux differentiability and 
elementary calculus, we get $\<\hat f\dual, h\>=0$. The assertions about the 
cases where $t$ is a bad point or a fan point are proved just as in 
Propositions 3.1 and 3.4.
\qed\enddemo

To finish this section we record an easy fact about uniform approximation. 
We may regard $\mu (f,u)$ as being the distance, calculated with respect to 
the norm $\norm{\cdot}$, from $f+f(t)\indic_{(t,u]}$ to the subspace $Z_u=
\{g\in \C_0(\Upsilon ):\supp g\subseteq (u,\infty)\}$, or as being the 
quotient norm $\norm{f+\indic_{(t,u]}+Z_u}$ in the quotient space 
$\C_0(\Upsilon )/Z_u$. Either way, it is very easy to establish the following 
lemma.

\proclaim{Lemma 3.5}
If $(f_n)$ is a sequence in $C_t$ which converges uniformly to $f$ then $\mu 
(f_n,\cdot)$ converges uniformly on $[t,\infty)$ to the limit $\mu(f,\cdot)$.
If $u$ is a bad point for all the functions $\mu (f_n,\cdot)$ then it is a bad 
point for  $\mu(f,\cdot)$. \endproclaim

\heading
4. LUR renormings
\endheading

The aim of this section is to establish a necessary and sufficient condition 
on the tree $\Upsilon $ for the existence of a locally uniformly convex norm 
on $\C_0(\Upsilon )$ equivalent to the supremum norm $\inorm{\cdot}$. In fact, 
this result can be deduced from the theorems in Sections 5 and 6 about 
strictly convex renormings and Kadec renormings, together with Troyanski's 
theorem. However, the direct proof that we are about 
to give is significantly shorter and easier to follow than such a roundabout 
approach, and we hope the reader will excuse the slight redundancy. We also 
show that the existence of a weakly LUR renorming of $\C_0(\Upsilon )$ implies 
the existence of a LUR renorming. It is an open problem whether wLUR 
renormability implies LUR renormability for arbitrary Asplund spaces, though 
this is the case for spaces which have a norm which is both
Fr\'echet-smooth and wLUR [\DGZ, Prop.~2.6]. 

\proclaim
{Theorem 4.1}
Let $\Upsilon $ be a tree. The following are equivalent:
\roster
\item there exists a locally uniformly convex renorming of $\C_0(\Upsilon )$;
\item there exists a weakly locally uniformly convex renorming of 
$\C_0(\Upsilon )$;
\item there exists an increasing function $\rho :\Upsilon \to\reals$ which is 
constant on no ever-branching subset of $\Upsilon $ and which has no bad 
points.
\endroster
\endproclaim

We already know from Propositions 3.2 and 3.4 that the function $\mu $ 
associated with an equivalent LUR norm on $\C_0(\Upsilon )$ can have no bad 
points and can be constant on no ever-branching subset. Thus the implication 
(1)$\implies$(3) is certainly true. We now pass to the implication 
(3)$\implies$(1). 

\proclaim
{Proposition 4.2}
Let $L$ be a locally compact space and let $U_i$ $(i\in I)$ be a family of 
open and closed subsets of $L$ such that, for each $i$, there is an equivalent 
LUR norm $\norm{\cdot}_i$ on $\C_0(U_i)$. For any bounded linear operator 
$S:\C_0(L)\to c_0(I)$, there is an equivalent norm $\norm{\cdot}$
on $\C_0(\Upsilon )$ which is locally uniformly convex at each $f$ such that 
the set $\{t\in L:f(t)\ne0\}$ is contained in the union $\bigcup\{U_i:i\in 
I\text{ and }(Sf)(i)\ne 0\}$. 
\endproclaim
\demo{Proof}
We may assume that $\norm{\cdot}_i\le \inorm{\cdot}$ for all $i$. For each 
non-empty finite subset $F$ of $L$ and each $f\in\C_0(L)$, we define
$$
\align
\phi (F;f) &= \biggl[\sum_{i\in F}((Sf)(i))^2\biggr]^{\frac12} ;\\
\psi (F;f) &= \biggl[\inorm{f.\indic_{L\less\bigcup_{i\in F}U_i}}^2+\sum_{i\in 
F}\norm{f\restriction U_i}_i^2\bigr]^{\frac12}.
\endalign
$$
Taking $I=I_m=\{F\subseteq I: \#F=m\}$ in Lemma 1.1, we obtain, for each 
positive integer $m$, a norm $\norm{\cdot}_m$ on $\C_0(L)$.  We set 
$\norm{f}^2=\sum_{m=1}^\infty 2^{-m} \norm{f}^2_m$ and claim that this new 
norm $\norm\cdot$ is locally uniformly convex at every $f$ satisfying the 
condition given above. 

Assume without loss of generality that $\inorm{f}=1$ and let $\epsilon $ be a 
positive real number. There is a finite subset $H_0$ of $I$ such that 
$(Sf)(t)\ne 0$ for all $i\in H_0$ and $\{s\in L: |f(s)|\ge \epsilon \} 
\subseteq \bigcup_{i\in H_0}U_i$. We may fix a positive real number $\alpha $ 
such that $H_0$ is contained in the finite set $H=\{t\in L:|(Sf)(t)|\ge 
\alpha\}$, noting 
that, for some $\eta >0$ we shall have $(Sf)(u)^2<\alpha ^2-\eta $ whenever 
$u\in I\less H$. Let $m=\#H$ and choose $k\in\nats$ so that $k>4\eta ^{-1}$. 
Notice that $\phi (H;f)=\sup_{F\in I_m}\phi (F;f)$ and that we have the 
following ``rigidity property'':
$$
G\in I_m\,,\ \phi (G;f)^2\ge \phi (H;f)^2-\eta \implies G=H.
$$

If $f_n$ are such that $2\norm{f_n}^2+2\norm{f}^2-\norm{f+f_n}^2\to 0$ then, 
we certainly have $2\norm{f_n}_m^2+2\norm{f}_m^2-\norm{f+f_n}_m^2\to 0$
and so, by Lemma~1.1, there exists a sequence $(F_n)$ in $I_m$ such that 
(among other things) $\phi (F_n;f)\to \sup_{F\in I_m}\phi (F;f)$ as $n\to 
\infty$. Because of the rigidity that we have just observed, $F_n$ must equal 
$H$ for all large enough $m$. Thus, by another of the conclusions of Lemma~1.1,
we see that 
$$ 
\psi 
(H;f)^2+\psi (H;f_n)^2-2\psi (H;\half(f+f_n))^2\to 0 $$ as $n\to\infty$, which 
in turn implies that
$$ \align 
\norm{(f_n)\restriction U_i}^2_i+\norm{(f)\restriction U_i}^2_i -
2\norm{\half(f+f_n)\restriction U_i}_i^2 &\to 0\\
\intertext{for all $i\in H$, and that} 
\inorm{f.\indic_{L\less\bigcup_{i\in 
H}U_i}}^2 + \inorm{f_n.\indic_{L\less\bigcup_{i\in H}U_i}}^2 -
2\inorm{\half(f+f_n).\indic_{L\less\bigcup_{i\in H}U_i}}^2 &\to 0. \endalign 
$$ 
The assumed local uniform convexity of the norms $\norm{\cdot}_i$ now gives 
us uniform convergence of $f_n$ to $f$ on the union $\bigcup_{i\in H}U_i$, 
while the second of the above limits certainly implies that 
$$ 
\inorm{f_n.\indic_{L\less\bigcup_{i\in H}U_i}}^2\to 
\inorm{f.\indic_{L\less\bigcup_{i\in H}U_i}}^2 
$$ 
Since $H$ contains our original $H_0$, we have
$\inorm{f.\indic_{L\less\bigcup_{i\in H}U_i}}<\epsilon $ and we thus see that 
$\inorm{f-f_n}\le 2\epsilon $ for all large enough $n$ \qed\enddemo 

\smallskip
We now define an operator $S$ which will here be used in an application of 
Proposition~4.2 and which will turn up again in later sections of the paper.
Recall from Section~3, that for an element of $\Upsilon $ which is a good point
for $\rho $ there exist 
$\delta _t>0$ and a finite (possibly empty) subset $F_t\subseteq t^+$ such that 
$$ 
\rho (u) \left\{\matrix\format \l & \l \\
         =\rho (t)  & \qquad{(u\in F_t)}\\ 
            \ge \rho (t)+\delta _t & \qquad{(u\in t^+\less F_t)}.\endmatrix
\right.
$$                            
For $f\in \C_0(\Upsilon )$ we define
$$
(Sf)(t)= \left\{\matrix\format \c & \l \\
\frac{\delta _t}{1+\#F_t}\biggl[f(t)-\sum_{u\in F_t} f(u)\biggr] & \quad 
\text{if $t$ is a good point} \\
0 & \quad \text{if $t$ is a bad point.}
\endmatrix\right.
$$          

Although for the immediate purposes of this section we are dealing with a tree 
equipped with a function $\rho $ with no bad points it is convenient for later 
applications to state the following two lemmas in greater generality.

\proclaim{Lemma 4.3}
Let $\Upsilon $ be a tree and let $\rho :\Upsilon \to \reals_+$ be a bounded 
increasing function.  The operator $S$ defined above is a bounded linear 
operator from $\C_0(\Upsilon )$ into $c_0(\Upsilon )$.
\endproclaim
\demo{Proof}
Certainly $S$ is a bounded linear operator from $\C_0(\Upsilon )$ into 
$\ell_\infty(\Upsilon )$. To show that $S$ takes values in $c_0(\Upsilon )$ it 
is enough to show that $S\indic_{(0,u]}\in c_0(\Upsilon )$ for all $u$. For 
$t\in (0,u]$ let us write $t'$ for the unique element of $t^+\cap (0,u]$; it 
is clear that $(S\indic_{(0,u]})(t)$ is non-zero if and only if $t=u$ or
$t\trl u$ and $t'\notin F_t$. Proceeding as in the proof of Proposition~2.7 
we calculate
$$
\align
\sum_{t\in\Upsilon }|(S\indic_{(0,u]})(t)|
&=\delta _u+\sum_{t\trl u,\ t'\ne F_t} \frac{\delta _t}{1+\#F_t}\\
&\le \delta_u+\sum_{t\trl u} (\rho (t')-\rho (t))\\
&\le \delta_u+\rho (u)
\endalign
$$
deducing that $S\indic_{(0,u]}$ is in $\ell_1(\Upsilon )$. 
\qed\enddemo

\proclaim{Lemma 4.4} Let $\Upsilon $ be a tree, let $\rho :\Upsilon \to \reals$ 
be a bounded increasing function which is continuous for the locally compact 
topology on $\Upsilon $, let $f$ be in $\C_0(\Upsilon )$ and let $s$ be an 
element of $\Upsilon $ with $f(s)\ne 0$. If the set $[s,\infty)\cap \rho ^{-
1}\rho (s)= \{t\in \Upsilon : t\trge s \text{ and } \rho (t)= \rho (s)\}$ has 
no ever-branching subset and contains only points that are good for $\rho $ 
then there exists $t$ in this set such that $(Sf)(t)\ne 0$.  
\endproclaim
\demo{Proof}
If no such $t$ exists we have $f(t)=\sum_{u\in F_t} f(u)$ 
whenever $t\trge s$ and $\rho (t)=\rho (s)$. Define a set $T$ by
$$
T=\{t\trge s: \text{ there exists } v\trge  t \text{ such that } \rho (v)=\rho 
(s) \text{ and }f(v)\ne 0 \}.
$$
By hypothesis, $T$ is not ever-branching so that there is some $t_0\in T$ such 
that $T\cap [t_0,\infty)$ is totally ordered. Now, for $t\in T\cap 
[t_0,\infty)$ we have  $f(t)=\sum_{u\in F_t} f(u)$ and at most one of the 
$f(u)$ can be non-zero. Thus $f$ is actually constant on $T\cap [t_0,\infty)$
and this constant value is non-zero because of the way we defined $T$. 
By the fact that $f\in \C_0(\Upsilon )$ and by the assumed continuity of $\rho 
$ we see that the set $T$, which can be written $T=\{t\trge t_0:\rho (t)=\rho 
(s) \text{ and } f(t)=f(t_0)\}$, is compact. So $T$ has a maximal point, $t$ 
say. This maximality implies that for any $v\trg t$ either $\rho (v)>\rho (t)$ 
or $f(v)=0$; in particular, either $F_t$ is empty or else $f(u)=0$ for all 
$u\in F_t$. In any case, $\sum_{u\in F_t} f(u)=0\ne f(t)$ and $(Sf)(t)\ne 0$ 
as claimed.\qed\enddemo

We can now establish the implication (3)$\implies$(1). We may assume that 
$\rho $ takes values in the real interval $(0,1)$.  We may furthermore suppose 
that $\rho $ is continuous for the locally compact topology of 
$\Upsilon $, that is to say that $\rho (t)=\sup _{s\trl t} \rho (s)$ whenever 
$t$ is a limit point of $\Upsilon $. This involves no loss of generality, 
since replacing $\rho (t)$ by $\lim_{s\uparrow t}\rho (s)$ at limit points 
$t$ does not introduce new bad points and does not create new ever-branching 
sets of constancy of $\rho $. For $t\in \Upsilon $ we define $U_t=(0,t]$, 
noting that $\C_0(U_t)$ has an equivalent LUR norm since $U_t$ is homeomorphic 
to the interval of ordinals $[0,r(t)]$. Lemma~4.3 shows that $S$ takes values 
in $c_0(\Upsilon )$.  Since all points of $\Upsilon $ are good and $\rho $ is 
constant on no ever-branching subset, Lemma~4.4 tells us that whenever 
$f(s)\ne 0$ there is some $t$ with $(Sf)(t)\ne 0$ and $s\in U_t$. Thus the 
hypotheses of Proposition 4.2 are satisfied and a LUR norm may be constructed 
on $\C_0(\Upsilon )$.

\smallskip
We now pass to the implication $\text{(2)}\implies\text{(3)}$. If we have a 
norm on $\C_0(\Upsilon )$ which is wLUR and not LUR then the associated 
function $\mu $ cannot be constant on an ever-branching subset (by 
Proposition 3.4) but it may quite easily have bad points. We show however that
it may be modified to give a new function $\rho $ without bad points.

\proclaim{Lemma 4.5}
Let $\norm{\cdot}$ be an equivalent weakly LUR norm on $\C_0(\Upsilon )$. The 
associated function $\mu :\Upsilon \to \reals$ has the following property:
whenever $(t_n)$ is an increasing sequence of points which are bad for $\mu $ 
and $t$ is another bad point with $t\ge \lim_n t_n$, we have $\mu (t)>\lim \mu 
(t_n)$.
\endproclaim
\demo{Proof}
Set $f_n=\indic_{(0,t_n]}$ and $f=\indic_{(0,t]}$. By Lemma 3.1 we have 
$\norm{f_n}=\mu (t_n)$ and $\norm{f}=\mu (t)$ and by definition of $\mu $
$\norm{\half(f+f_n)} \ge \mu (t_n)$. If it were the case that $\mu 
(t)=\lim_n\mu (t_n)$ then we should have  $\norm{f_n}\to \norm{f}$ and 
$\norm{\half(f+f_n)}\to \norm{f}$, which contradicts the wLUR property, since 
$f(t)=1, $ and $f_n(t)=0$ for all $n$.
\qed\enddemo

\smallskip
We can now complete the proof of the theorem. We already know that the 
function $\mu $ associated with a wLUR norm on $\C_0(\Upsilon )$ 
is strictly increasing on the set of its 
bad points (Proposition 3.3) and that it satisfies the condition established in 
the above lemma. Together these two facts tell us that for any point $t$ 
which is bad for $\mu $ we have 
$$
\mu (t)>\sup\{\mu (s): s \text{ is bad for $\mu $ and } s\trl t\}.
$$
We define
$$
\rho (t)=\mu(t) + \sup\{\mu(s): s \text{ is bad for $\mu $ and } s\trl t\},
$$
and note first that a point which is good $\mu$ must certainly be good for 
$\rho $. If $t$ is bad for $\mu $, on the other hand, then for any $u\in t^+$ 
we have 
$$
\rho (u)=\mu (u)+\mu (t) \ge 2\mu (t)>\rho (t)
$$
by what we have just observed. Thus $t$ is good for $\rho $. The fact that 
$\rho $ is constant on no ever-branching subset follows immediately from the 
same property for $\mu $, and so the proof of Theorem 4.1 is complete.

\heading
5. Strictly convex renormings
\endheading

In this section we establish a necessary and sufficient condition for the 
existence on $\C_0(\Upsilon )$ of an equivalent strictly convex norm and we 
further show that this condition yields the existence of an equivalent norm 
that is midpoint-locally uniformly convex (MLUR). The 
condition may be regarded as a mild weakening of $\reals$-embeddability and 
the theorem implies the result about $\reals$-embeddable trees that we 
presented as Proposition~2.7.

\proclaim{Theorem 5.1}
For a tree $\Upsilon $ the following are equivalent: 
\item{\rm (1)} $\C_0(\Upsilon )$ admits an equivalent strictly convex norm;
\item{\rm (2)} $\C_0(\Upsilon )$ admits an equivalent MLUR norm;
\item{\rm (3)} there exists on $\Upsilon $ an 
increasing real-valued function $\rho $ satisfying the two conditions:
\itemitem{\rm (i)} $\rho$ is constant on no ever-branching subset of $\Upsilon$;
\itemitem{\rm (ii)} for any $s\in \Upsilon $ there is at most one bad point 
$u$ with $s\trle u$ and $\rho (u)=\rho (s)$.
\item{\rm (4)} there is a bounded linear injection from $\C_0(\Upsilon )$ into 
some space $c_0(I)$.
\endproclaim

The fact that (1) implies (3) has already been established in Propositions 1.8 
and 2.2, and certainly either of (2) or (4) implies (1). We have to show that 
(3) implies (2) and (4). So let $\rho :\Upsilon \to (0,1)$ be an increasing 
function sastisfying (i) and (ii). 

\smallskip
To establish (4) we define $R:\C_0(\Upsilon )\to c_0(\Upsilon )$, 
$S:\C_0(\Upsilon )\to c_0(\Upsilon )$ by
$$
\align
(Rf)(t) &= \left\{\matrix 
                   \format \c & \l \\  
\rho(t)-\rho(t^-))f(t) &\qquad\qquad\qquad\qquad \text{if $t\in \Upsilon ^+$}\\
            0 &\qquad\qquad\qquad\qquad \text{otherwise} \endmatrix\right.\\
\ &\ \\
(Sf)(t) &= \left\{\matrix\format \c & \l \\
 \delta_t (\#F_t)^{-1} \bigl(f(t)-\sum_{u\in F_t}f(u)\bigr) 
              &    \qquad \text{if } t \text{ is good,}\\
                  0 & \qquad \text{otherwise.}
          \endmatrix\right.
\endalign
$$                                                 
We shall show that the operator $R\oplus S:\C_0(\Upsilon )\to c_0(\Upsilon) 
\oplus c_0(\Upsilon )$ is injective.

Suppose then  that $f\in \C_0(\Upsilon )$ is such that $Rf=Sf=0$. If $f\ne 0$ 
we choose a minimal element $t$ of $\Upsilon $ with $f(t)\ne 0$. We notice 
straightaway that by Lemma 4.4 there must be at least one bad point $u$ with 
$u\trge t$ and $\rho (u)=\rho (t)$. By continuity of $f$ and minimality of 
$t$, it cannot be that $t$ is a limit element and so $t$ has an immediate 
predecessor $s=t^-$; since $Rf=0$ it must be that $\rho (t)=\rho(s)$. By 
minimality of $t$ we have $f(s)=0$ and since $Sf=0$ we have $\sum_{v\in 
F_{s}\less\{t\}}f(v)=-f(t)\ne0$. Let $t'$ be any element of $F_{s}\less\{t\}$ 
with $f(t')\ne 0$. By Lemma 4.4 applied this time to $t'$ there exists at 
least one bad point $u'$ with $u'\trge t'$ and $\rho (u')=\rho (t')=\rho (s)$.  
The existence of two bad points $u,u'$ with $u,u'\trg s$ and $\rho (u)=\rho 
(u')=\rho (s)$ contradicts property (ii) of the function $\rho$. 

\smallskip
We now set about constructing a MLUR renorming of $\C_0(\Upsilon )$. To get an 
idea of how one can have such a renorming when no LUR renorming exists one can 
start with the observation that in the unit sphere of the space $\C[0,1]$, 
under the supremum norm, there are no points of local uniform convexity, 
although the constant functions $\pm1$ are points of mid-point local uniform 
convexity. The following lemma presents a small generalization of this remark.
For a bounded function $f$ defined on a set $V$, we write 
$\text{osc}(f)=\sup_{s,t\in V}|f(s)-f(t)|$ and we define the oscillation norm 
$\norm{\cdot}_{\text{osc}}$ by
$$
\norm{f}_{\text{osc}}^2 = \inorm{f}^2 + \text{osc}( f)^2.
$$

\proclaim{Lemma 5.2}
An element $g$ of $\ell_\infty(V)$ is a point of mid-point local uniform 
convexity of $\norm{\cdot}_{\text{osc}}$ if and only if $g$ takes at most two 
values. If $g$ is ``within $\epsilon $ of being two-valued'', in the sense 
that there exist real numbers $a,b$ such that $g$ takes values in $[a-\epsilon 
,a+\epsilon ]\cup[b-\epsilon ,b+\epsilon ]$, then $g$ has the following 
property:
$$
\norm{g+h}_{\text{osc}}^2+\norm{g-h}_{\text{osc}}^2-
2\norm{g}_{\text{osc}}^2<\epsilon ^2 \implies \norm{h}<4\epsilon .
$$
\endproclaim
\demo{Proof}
We give a proof only of the second assertion. It is not hard to see that the 
hypothesis about $g\pm h$ implies that
$$
\align 
\inorm{g+h}^2+\inorm{g-h}^2-2\inorm{g}^2 &<  \epsilon ^2,\\
\text{osc}(g+h)^2+\text{osc}(g-h)^2-2\text{osc}(g)^2 
&<  \epsilon ^2,
\endalign
$$
and that these inequalities imply 
$$
\inorm{g\pm h}< \inorm{g}+\epsilon ,\quad \text{osc}(g\pm h) 
< \text{osc}(g)+\epsilon .
$$
Now let $A=\sup g[V]$ and $B=\inf g[V]$ and assume without loss of generality 
that $\inorm{g}=A$, $\osc{g}=A-B$. By what we have just proved we have 
$$
\inorm{g\pm h}< A+\epsilon ,\quad \text{osc}(g\pm h) 
< A-B+\epsilon .
$$
If $s\in V$ is such that $g(s)\ge A-2 \epsilon $ it is now clear that 
$|h(s)| $ can be at most $3\epsilon $. It is clear also that if $s$ is such 
that $g(s)$ is suficiently close to $A$ then $|h(s)|<\epsilon $. 
We now choose such an $s$ and 
consider  $t$  such that $g(t)\le B+2\epsilon $. We have $|h(t)|=|g(s)-(g(t) 
+\sigma  h(t))|-|g(t)-g(s)|$ for a suitable choice of sign $\sigma =\pm1$. 
Thus 
$$ 
\align
|h(t)| &\le \epsilon + |(g+\sigma h)(s)-(g+\sigma h)(t)|-|g(s)-g(t)| \\
       &\le \epsilon + \text{osc}(g+\sigma h) - (A-B)+2\epsilon \le 4\epsilon.
\endalign
$$ 
We have now shown that $|h|<3\epsilon $ on the set where $A-2\epsilon \le g\le 
A$ and that $|h|<4\epsilon $ on the set where $B\le g\le B+2\epsilon $. The
hypothesis that $g$ is almost two-valued tells us these are the only 
possibilities and hence that $\inorm h\le 4\epsilon $, as claimed.
\qed\enddemo

We now present a proposition which is analogous to Proposition~4.2.

\proclaim
{Proposition 5.3}
Let $L$ be a locally compact space, let $(U_i)_{i\in I}$,
$(V_i)_{i\in I}$ be  families of 
open and closed subsets of $L$ and let $T:\C_0(\Upsilon )\to c_0(I)$ be a 
bounded linear operator. Assume that, for each $i$, there is an equivalent 
MLUR norm $\norm{\cdot}_i$ on $\C_0(U_i)$. Assume further that, 
for all $f\in \C_0(L)$ and all $\epsilon >0$, the 
set $\{t\in L:f(t)\ne0\}$ is contained in the union 
$$
\align
\bigcup\{U_i:(Tf)(i)\ne 0\}
\cup\bigcup\{V_i: &(Tf)(i)\ne 0\\
&\text{\sl\ and } f\restriction V_i \text{ is 
within $\epsilon $ of being two-valued}\} .
\endalign
$$ 
Then $\C_0(L)$ admits a MLUR renorming.
\endproclaim
\demo{Proof}
For technical reasons it is convenient to introduce a total order on the index 
set $I$. When $m$ is a positive integer and $F$ is a subset of $I$ with 
$\#F=m$ there is thus a fixed increasing order $i_0<i_1<\cdots<i_{m-1}$ in 
which to write the elements of $F$. If $\pi \in \{0,1\}^m$ we write $F^\pi $ 
for the subset $\{i_j:j<m \text{ and } \pi (j)=1\}$ of $F$. 

Now, for each positive integer $m$, each subset $F$ of $I$ with $\# F=m$ 
and each $\pi \in \{0,1\}^m$, we define
$$
\align
\phi (F;f) &= \biggl[\sum_{i\in F}((Tf)(i))^2\biggr]^{\frac12} ;\\
\psi (\pi, F;f) &= \biggl[\inorm{f.\indic_{L\less(\bigcup_{i\in F}U_i
\cup\bigcup_{i\in F^\pi }V_i}}^2+\sum_{i\in 
F}\norm{f\restriction U_i}_i^2
+\sum_{i\in F^\pi}\norm{f\restriction V_i}_{\text{osc}}^2\bigr]^{\frac12}.
\endalign
$$
We now apply the norm construction of Lemma~1.1 obtaining a norm $\norm\cdot$ 
such that
$$
\align
\norm{f}^2 &= \inorm{f}^2 + \sum_{m=1}^\infty\infty2^{-m-2^m}
\sum_{\pi \in \{0,1\}^m}\theta (m,\pi ;f),\\
\intertext{where}
\theta (m,\pi ;f) &= \sum_{l=1}^\infty 2^{-l}
\sup{\#F=m}
\biggl[\phi (F;f)^2+2^{-l}\psi (\pi ,F,f)^2\biggr].
\endalign
$$

For any $f\in \C_0(L)$ and any $\epsilon >0$ there exists a finite 
$H_0\subset I$ such that $\{t\in L:|f(t)|\ge \epsilon \}$ is contained in
the union of the $U_i$ $(i\in H_0)$, together with those $V_i$ $(i\in H_0)$
for which $f\restriction V_i$ is within $\epsilon $ of being two-valued. As in 
the proof of Proposition~4.2, we find $\alpha >0$ such that $H_0\subseteq H= 
\{i\in I: |(Tf)(i)|\ge \alpha \}$ and we set $m=\#H$. As a small refinement, 
we let $\pi \in \{0,1\}^m$ be the function such that $H^\pi $ contains exactly 
those $i\in H$ such that $f\restriction V_i$ is within $\epsilon $ of being 
two-valued. 

Now let $(h_n)$ be a sequence in $\C_0(L)$ such that
$\norm{f+h_n}^2+\norm{f-h_n}^2-2\norm{f}^2\to0$. For the $m$ and $\pi $ which 
we fixed above, we have
$$
\theta (m,\pi ;f_n)+\theta (m,\pi ;f)-2\theta (m,\pi ;\half(f+f_n))\to 0
$$
as well. If $(F_n)$ is the sequence given by the conclusion of Lemma~1.1 then, 
as in the proof of Proposition~4.2, we must have $F_n=H$ for all large enough 
$n$. We also have
$$
\psi (H,\pi ;f+h_n)^2+\psi (H,\pi ;f+h_n)^2-2\psi (H,\pi ;f)\to 0
$$
as $n\to \infty$, which implies in turn that
$$
\align
\norm{(f+h_n)\restriction U_i}^2_i+\norm{(f-h_n)\restriction U_i}^2_i
-2\norm{f\restriction U_i}_i^2 &\to 0\quad (i\in H),\\
\norm{(f+h_n)\restriction V_i}^2_{\text{osc}}
+\norm{(f-h_n)\restriction V_i}^2_{\text{osc}}
-2\norm{f\restriction V_i}_{\text{osc}}^2 &\to 0\quad (i\in H^\pi ),\\
\inorm{(f\pm h_n).
         \indic_{L\less(\bigcup_{i\in H}U_i\cup\bigcup_{i\in H^\pi }V_i)}}\to
\inorm{f.\indic_{L\less(\bigcup_{i\in H}U_i\cup\bigcup_{i\in H^\pi }V_i)}}&\le 
\epsilon.
\endalign
$$
These limits give us, for all large enough $n$,
$$
\align
|h_n|<\epsilon  &\text{ on $U_i$ $(i\in H)$ by the MLUR property of 
$\norm{\cdot}_i$,}\\
|h_n|<4\epsilon &\text{ on $V_i$ $(i\in H^\pi) $ by Lemma~5.2,}\\
|h_n|<2\epsilon &\text{ on the rest of $L$}.
\endalign
$$
We have proved that $\inorm{h_n}\to 0$ and thus that $\norm{\cdot}$ is MLUR.
\qed\enddemo

\smallskip
We are now ready to complete the proof of Theorem~5.1 by showing that, subject 
to condition (3), we can find $I,T,U_i\,,V_i$ to which Proposition~5.3 may be 
applied. 

We take $T$ to be the operator $R\oplus S:\C_0(\Upsilon )\to c_0(\Upsilon 
)\oplus c_0(\Upsilon )$ that we have already looked at in this section. For 
$t\in \Upsilon $ we define $a(t)$ to be the minimal element of $\{s\in (0,t]:
\rho (s)=\rho (t)\}$ and, if there is a bad point $u$ with $u\trge a(t)$ and 
$\rho (u)=\rho (t)$ we define $b(t)$ to be this point (which is unique because 
of condition (ii)). We set $U_t=(0,t]\cup(a(t),b(t)]$ and $V_t=[t,\infty)$, 
noting that since $U_t$ is homeomorphic to an interval of ordinals there is 
certainly a LUR renorming of $\C_0(U_t)$. 

\proclaim{Lemma 5.4}
Let $\Upsilon $ and $\rho $ satisfy conditions (i) and (ii) and suppose 
further that $\rho $ is continuous for the locally compact topology of 
$\Upsilon $. Let $f$ be in $C_0(\Upsilon )$ and let $s$ be an element of 
$\Upsilon $, with $f(s)\ne 0$,  which is not in the union 
$\bigcup\{U_t:(Rf)(t)\ne 0 \text{ or } (Sf)(t)\ne 0\}$. Then the restriction 
of $f$ to $[a(s),\infty)$ is a constant multiple of the indicator function
$\indic_{[a(s),b(s)]}$. Moreover, $a(s)$ is a limit element of $\Upsilon $ and 
there is a sequence $(r_n)$ increasing to $a(s)$ in $\Upsilon $ such that
$(Rf)(r_n)\ne 0$ for all $n$.
\endproclaim
\demo{Proof}
We start by noting that certainly there is no $t\trge s$ with $(Sf)(t)\ne 0$.
It follows from Lemma~4.4 that $[s,\infty)\cap \rho ^{-1}\rho (s)$ cannot 
contain only good points, so it must be that $b(s)$ exists and that 
$s\in [a(s),b(s)]$. Now for any $t$ with $t\trge a(s)$ and $\rho (t)=\rho (s)$ 
we have $a(t)=a(s)$ and $b(t)=b(s)$, so that $s\in U_t$ whence, by hypothesis, 
$(Sf)(t)=0$. Using Lemma~4.4 again, we can now see that on the set 
$[a(s),\infty)\cap \rho ^{-1}\rho (s)$ our function $f$ can be nonzero only on 
$[a(s),b(s)]$, and since $(Sf)(t)=0$ for all $t\in [a(s),b(s))$ it must be 
that $f$ is constant on $[a(s),b(s)]$. 

If $a(s)$ were not a limit element of $\Upsilon $ we should have $\rho 
(a(s))>\rho (a(s)^-)$ by the minimality in the definition of $a(s)$ and since 
$f(a(s))=f(s)\ne 0$ we should have $(Rf)(a(s))\ne 0$ and  
$s\in [a(s),b(s)]\subseteq U_{a(s)}$, contrary to assumption. Since $a(s)$ is 
a limit point with $\rho (a(s))>\rho (r)$ when $r \trl a(s)$ there is a 
sequence $(r_n)$ increasing to $a(s)$ consisting of successor elements with 
$\rho (r_n)>\rho (r_n^-)$. Because $f$ is continous it is non-zero on some 
neighbourhood of $a(s)$ and so we may assume that $f(r_n)\ne 0$ for all $n$.

Finally we have to show that if $s'\trg a(s)$ and $\rho (s')>\rho (s)$ then 
$f(s')\ne 0$. We note that for such an $s'$ we have $s\in U_t$ whenever $s'\in 
U_t$ so that $(Sf)(t)=0$ whenever $s'\in U_t$. Hence, if $f(s')$ is not 0, 
all of the above analysis is applicable to $s'$ and in particular $a(s')$ is a 
limit element, with $f(s')\ne 0$. Also there exists a sequence $(r'_n)$ 
increasing to $a(s')$ with $(Rf)(r_n')\ne 0$ for all $n$. But this cannot be, 
since we have $s\in (0,r_n']\subseteq U_{r_n'}$ for suitably large $n$.
\qed\enddemo

The lemma we have just proved really does all the work for us since, given 
$s$ with $f(s)\ne 0$ and $s\notin \bigcup\{U_t: (Rf)(t)\ne 0 \text{ or } 
(Sf)(t)\ne 0\}$, we get a sequence $(r_n)$ increasing to $a(s)$ and consisting 
of points with $(Rf)(r_n)\ne 0$. Given $\epsilon >0$ we may find $n$ such that
$|f(r)-f(a(s))|<\epsilon $ for all $r\in [r_n,a(s)]$ and such that 
$|f(t)|<\epsilon $ when $t\in [r_n,\infty)\less [a(s),\infty)$, the second of 
these properties resulting from the compactness of the set where $|f|\ge 
\epsilon $. Together with what we have already proved about $f$ restricted to 
$[a(s),\infty)$, this shows that on $V_{r_n}=[r_n,\infty)$ the function $f$ is 
indeed within $\epsilon $ of being two-valued.

We have now shown that the hypotheses of Proposition~5.3 are satisfied and we 
have thus completed the proof of the theorem.

\heading 
6. Kadec renormings 
\endheading

\proclaim{Theorem 6.1} Let $\Upsilon$ be a tree. The Banach space 
$\C_0(\Upsilon )$ admits an equivalent norm with the Kadec Property if and 
only if there exists an increasing function $\rho :\Upsilon \to \reals$ with 
no bad points. \endproclaim 

We  have already seen that if $\C_0(\Upsilon )$ admits a Kadec norm then the 
associated $\mu $ function has no bad points. We shall devote this section to 
the construction of a Kadec renorming, starting from a function $\rho $ with 
no bad points. As in the proof of Theorem~4.1, may suppose that $\rho $ takes 
values in $(0,1)$ and is continuous for the locally compact topology of 
$\Upsilon $. 

For $r\in\Upsilon $, $f\in\C_0(\Upsilon )$ and $1\le l\in\nats$ we make the 
following definitions:
$$
\align
\Delta ^\pm(f;r) = &\inf \{\inorm{g\mp f\restriction_{[r,\infty)}}: g\in 
\C_0[r,\infty) \text { and } g \text{ is decreasing}\};\\
\Alpha_l(f,r) = \frac1l&\sup\biggl\{\sum_{k=1}^l |f(s_k)|:\{s_1,s_2,\dots,s_l\} 
\text{ is an antichain in } [r,\infty)\biggr\}.
\endalign
$$
Notice that $\Alpha_1(f,r) = \inorm{f\restriction_{[r,\infty)}}$.
We recall our standard notation for good points:
$$
F_t=\{u\in t^+:\rho (u)=\rho (t)\}, \qquad \delta _t = \inf\{\rho (v)-\rho 
(t): v\in t^+\less F_t\}.
$$

\proclaim{Lemma 6.2}
There are uniquely determined functions $\Phi ,\Psi ,\Theta ,\Omega,\Sigma  
:\C_0(\Upsilon )\times\bigl(\Upsilon \cup\{0\}\bigr)\to \reals^+$ and $\Xi :
\C_0(\Upsilon )\times\bigl(\Upsilon \cup\{0\}\bigr)\times
\bigl(\Upsilon \cup\{0\}\bigr)\to \reals^+$ 
which satisfy the inequalities
$$
\Phi(f;s) ,\Psi(f;s) ,\Theta (f;s),\Omega(f;s),\Sigma (f;s), \Xi (s,t) \le 
\inorm{f}
$$ 
for $f\in \C_0(\Upsilon )$ and $s\trle t\in \Upsilon $ as well as 
the identities:
$$
\align
7\Phi (f;s) 
&= \Delta ^+(f;s) +\Delta ^-(f;s)+\sum_{l=1}^\infty 2^{-l} \Alpha_l(f;s)\\ 
&\qquad\qquad\qquad\qquad\qquad+\Sigma (f;s)+\Psi(f;s) +\Theta (f;s)+\Omega(f;s);\\
2\Sigma (f;s) &= \sum_{m,l=1}^\infty \!2^{-m-l}\sup_{s\trle 
t}\biggl[|f(t)|+2^{-m}\Xi (f;s,t)+2^{-l}\Phi (f,t)\biggr];\\
4\Psi (f;s) &=  \sum_{m,l=1}^\infty\!2^{-m-l}
\!\!\!\!\!\!\!
\sup_{s\trle t\trle u\trle\infty } \biggl[\,\bigl|f(s)-f(t)+f(u)\bigr| 
+2^{-m}\Xi (f;s,t)+2^{-l}\Phi (f;t)\biggr] ;\\
2\Theta (f;s) &= \sum_{m,l=1}^\infty\!2^{-m-l}
\sup_{s\trle t} \, \frac{1}{1+\#F_t}\biggl[ 
\delta _t\Bigl|f(t)-\!\!\sum_{u\in F_t}f(u)\Bigr|\\
&\qquad\qquad\qquad\qquad\qquad\qquad
+\sum_{u\in F_t}\Bigl(2^{-m}\Xi (f;s,u)+2^{-l}\Phi (f;u)\Bigr)\biggr] ;\\
2\Omega (f;s) &= \sum_{m,l=1}^\infty\!2^{-m-l}
\!\!\!\sup_{s\trl t\in\Upsilon ^+}
\biggl[ (\rho (t)-\rho (t^-))|f(t)|+2^{-m}\Xi (f;s,t)
+2^{-l}\Phi (f;t)\biggr];\\
2\Xi (f;s,t)&= \onorm{f\restriction_{(s,t]}}+\Phi 
(f.\indic_{\Upsilon \less(0,t,\infty)};r).
\endalign
$$
As functions of $f$, all of $\Phi ,\Sigma ,\Psi ,\Theta ,\Omega ,\Xi $ are 
$\tau_{\text p}$-lower semicontinuous seminorms on $\C_0(\Upsilon )$. They are 
continuous for the locally compact topology of $\Upsilon $ as functions of
their arguments $s,t$. All of these functions are decreasing in $s$, and $\Xi 
$ is increasing in $t$.
\endproclaim
\demo{Proof}
The existence and uniqueness are immediate consequences of Banach's 
fixed-point theorem in a suitable space of functions. The other assertions are 
straightforward to verify. 
\qed
\enddemo

\smallskip
It is clear that $\norm{f}=\Phi (f;0)$ defines a norm on $\C_0(\Upsilon )$ 
such that $\quarter \inorm{f}\le \norm{f}\le \inorm{f}$.  We shall show that 
$\norm{\cdot}$ has the Kadec Property.

\proclaim{Lemma 6.3} Let $f_n$ and $f$ in $\C_0(\Upsilon )$ and $r\in \Upsilon $ 
be such that $f_n\to f$ pointwise and $\Phi (f_n;r)\to \Phi (f;r)$. Then
$\Sigma (f_n;r),\ \Psi (f_n;r),\ \Theta (f_n;r),\ \Omega (f_n;r),\ \Delta 
^\pm(f_n;r),$ and $\Alpha_l(f_n;r)$ $(l\ge 1)$ converge as $n\to \infty$ to the 
limits $\Sigma (f;r)$, $\Psi (f;r)$, $\Theta (f;r)$, $\Omega (f;r)$, 
$\Delta ^\pm(f;r),$ and $\Alpha_l(f;r)$, respectively.  In particular,
$\inorm{f_n\restriction_{[r,\infty)}}\to\inorm{f\restriction_{[r,\infty)}} $ .
\endproclaim
\demo{Proof}
This follows immediately from the $\tau_{\text p}$-lower semicontinuity of the 
terms in the sum defining $\Phi $.\qed
\enddemo

\proclaim{Lemma 6.4} Let $f_n$ and $f$ in $\C_0(\Upsilon )$ and $r\in \Upsilon $ 
be such that $f_n\to f$ pointwise and $\Phi (f_n;r)\to \Phi (f;r)$. Let 
$\epsilon $ be any real number with $0<\epsilon \le 
\inorm{f\restriction_{[r,\infty)}}$. There exists $s\in [r,\infty)$ such that 
\roster          
\item \ $|f(s)|=\inorm{f\restriction_{[s,\infty)}}\ge \epsilon $ ;
\item \ $f_{n}\to f$ uniformly on $(r,s]$;
\item \ $\Phi (f_{n}\indic_{\Upsilon \less(0,s,\infty)};r)
     \to \Phi (f\indic_{\Upsilon \less(0,s,\infty)};r)$;
\item \ $\Phi (f_{n};s)\to \Phi (f;s)$.
\endroster 
\endproclaim 
\demo{Proof}
We may assume that $\inorm{f}<1$, so that the same is true for any of the 
functions $\Phi ,\Sigma $ and so on, associated with $f$.
The set $M$ of maximal elements of $\{t\in [r,\infty): |f(t)|\ge\epsilon \}$ 
is finite and we can choose a natural number $m$ such that $|f(s)|< \epsilon -
2^{-m}$ whenever $s\in [r,\infty)\less \bigcup_{t\in M} [r,t]$. 

By Lemma~6.3, $\Sigma (f_n;r)\to \Sigma (f;r)$ as $n\to \infty$ and it 
follows similarly that 
$$
 \sum_{l=1}^\infty \!2^{-l}\sup_{r\trle 
s}\biggl[|f_n(s)|+2^{-m}\Xi (f_n;r,s)+2^{-l}\Phi (f_n,s)\biggr]
$$
tends to 
$$
 \sum_{l=1}^\infty \!2^{-l}\sup_{r\trle 
s}\biggl[|f(s)|+2^{-m}\Xi (f;r,s)+2^{-l}\Phi (f,s)\biggr].
$$
Now by applying Lemma~1.2 with $I=[r,\infty)$, $\phi _s(g)= |g(s)|+2^{-m}\Xi 
(g;r,s)$, $\psi _s(g)=\Phi (g;s)$ we see that there  is a sequence $(s_n)$ in 
$[r,\infty)$ such that $|f(s_n)|+2^{-m}\Xi (f;r,s_n)$ and $|f_n(s_n)|+2^{-
m}\Xi (f_n;r,s_n)$ and $\sup_{t\in [r,\infty)}|f_n(t)|+2^{-
m}\Xi (f_n;r,t)$ all converge as $n\to \infty$ to the limit $\sup_{t\in 
[r,\infty)}\Bigl[|f(t)|+2^{-m}\Xi (f;r,t)\Bigr]$, and moreover such that $\Phi 
(f_n;s_n)-\Phi (f;s_n)$ tends to zero. We may further suppose that  $(s_n)$
converges in the sequentially compact space $\Upsilon \cup\{\infty\}$. Its 
limit must be a point $s$ of $\Upsilon $, since $f(s_n)$ does not tend to 
zero as $n\to \infty$.  In fact, $s$ is a point at which 
the supremum $\sup_{t\in [r,\infty)}\Bigl[|f(t)+2^{-m}\Xi (f;r,t)\Bigr]$ is 
attained and hence, by our original choice of  $m$, it must be in 
$\bigcup_{t\in M}[r,t]$. If $t\in M$ is such that $s\trle t$ then we have
$|f(s)|\ge |f(t)|+\Xi (f;r,t)-\Xi (f;r,s)\ge |f(t)|$, since $\Xi (f;r,\cdot)$ 
is increasing.  Thus $|f(s)|\ge \epsilon $. 

It is easy to see that (4) holds. Indeed, we have $\liminf \Phi (f_{n};s)\ge 
\Phi (f;s)$ by $\tau_{\text{p}}$-lower semicontinuity. On the other hand, 
since each $\Phi(g;t)$ is a continuous decreasing function of $t$, we have 
$\Phi (f;s_{n})\to \Phi (f;s)$ and $\Phi (f_{n};s_{n})\ge \Phi (f_{n};s)$. 
Using these facts and remembering that $\Phi (f_{n};s_{n})-\Phi (f;s_{n})$ 
tends to zero, we see that $\limsup\Phi (f_{n};s)\le \limsup\Phi 
(f_{n};s_{n})=\lim \Phi (f;s_{n})=\Phi (f;s)$. 

We next show that $\Xi (f_{n};r,s)\to \Xi (f;r,s)$ as $n\to\infty$. By 
$\tp$-lower semicontinuity, we have $\liminf\Xi (f_{n};r,s)\ge \Xi (f;r,s) $.
On the other hand, 
$$
\align
\limsup \Xi (f_{n};r,s) &\le \limsup\Bigl\{\sup_{t\in[r,\infty)}
\Bigl[2^m|f_{n}(t)|+ \Xi (f_{n};r,t)\Bigr] - 2^m|f_{n}(s)|\Bigr\}\\
    &= \limsup\Bigl\{2^m|f(s_{n})|+ \Xi (f;r,s_{n})\Bigr\} - 2^m|f(s)|\\
    &= 2^m|f(s)|+ \Xi (f;r,s) - 2^m|f(s)|= \Xi (f;r,s).
\endalign
$$
The definition of $\Xi $ and the usual lower semicontinuity argument now gives 
us the condergence of $\Phi (f_{n}\indic_{\Upsilon \less(0,s,\infty)};r)$ to 
$\Phi (f\indic_{\Upsilon \less(0,s,\infty)};r)$, which is (3), as 
well as the convergence of $ \onorm{f_n\restriction_{(r,s]}}$ to  
$\onorm{f\restriction_{(r,s]}}$, from which (2) follows because 
$\onorm{\cdot}$ is a Kadec norm. \qed \enddemo

The theorem will follow from the following proposition, the form of which is 
chosen to simplify a proof by induction. 

\proclaim{Proposition 6.5}
Let $\epsilon$ be a positive real number, let $q$ be in $\Upsilon\cup\{0\} $ 
and let $f$ be an element of $\C_0(\Upsilon )$.  If the sequence $(f_n)$ in 
$\C_0(\Upsilon )$  tends pointwise to $f$ and is such that 
$\Phi (f_n;q)\to \Phi (f,q)$, then 
$\limsup_{n\to\infty}\inorm{(f-f_n)\restriction_{[q,\infty)}}< 2\epsilon $. 
\endproclaim

We define 
$$
m(g;r;\epsilon )= \max\{\#A:A \text{ is an antichain in $[r,\infty)$ and 
$|g(t)|\ge  \epsilon $ for all } t\in A\},
$$
and note that if $m(f;q;\epsilon)=0$ then $\inorm{f\restriction_{[q,\infty)}}< 
\epsilon $. Since $\Phi (f_n;q)\to \Phi (f,q)$, we have 
$\inorm{f_n\restriction_{[q,\infty)}}=\Alpha_1(f_n;q)\to \Alpha_1(f;q)= 
\inorm{f\restriction_{[q,\infty)}}$, so that
$\inorm{f_n\restriction_{[q,\infty)}}$ is smaller than $\epsilon $ for all 
sufficiently large $n$. 

We now suppose that $m(f,q,\epsilon )>0$ and assume inductively that our 
result is true for any pair $g,r$ such that $m(g;r;\epsilon )< 
m(f;q;\epsilon)$.  

We write $K$ for the subset
$$
\align
K= \Bigl\{r\in [q,\infty):& \inorm{f\restriction_{[r,\infty)}}\ge 
\epsilon \text{ and } \lim_{n\to\infty}\Phi (f_{n},r)= \Phi (f,r) \\
&\text{ and }\limsup_{n\to\infty}
\inorm{(f-f_{n})\restriction_{[q,\infty)\less[r,\infty)}} 
<2\epsilon \Bigr\} .
\endalign 
$$ 
and set
$$
H=\{s\in K: |f(s)|=\inorm{f\restriction_{[s,\infty)}}\}.
$$

Evidently $K$ is non-empty since $q\in K$. Our first task is to show that $H$ 
is non-empty.  

\proclaim{Claim 1}
For every $r\in K$ there exists $s\in H\cap[r,\infty)$.
\endproclaim
\demo{Proof}
Since $\Phi (f_n;r)\to \Phi (f;r)$  we can apply Lemma~6.4 (with the 
``$\epsilon $'' of that lemma equal to $\inorm{f\restriction_{[r,\infty)}}$)
getting an element $s$ of $[r,\infty)$. We have 
$|f(s)|=\inorm{f\restriction_{[r,\infty)}}\ge \epsilon $ and $\Phi 
(f_{n};s)\to \Phi (f;s)$, by (1) and (4) in that Lemma. Now (3) tells us that 
as $n\to\infty$ the quantity $\Phi (f_{n}.\indic_{\Upsilon 
\less(0,s,\infty)};r) $ tends to the limit $\Phi (f.\indic_{\Upsilon 
\less(0,s,\infty)};r)$. Since $|f(s)|\ge \epsilon $ and $\Upsilon 
\less(0,s,\infty)$ contains no points comparable with $s$ we have 
$m(f.\indic_{\Upsilon \less(0,s,\infty)},r,\epsilon )<m(f,q,\epsilon)$ so that 
our inductive hypothesis yields
$\limsup\inorm{(f-f_n)\restriction _{[r,\infty)\less [r,s,\infty)}}<2\epsilon$. 
Finally, the assertion (2) of Lemma~6.4 tells 
us that $f_{n}\to f$ uniformly on $[r,s]$, so that certainly
$\limsup\inorm{(f-f_n)\restriction _{[r,s]}}<2\epsilon $. \qed\enddemo 

We may suppose that $H$ is  totally ordered by $\trle $. Indeed, if $r,s\in 
H$ are incomparable, we have $[q,\infty)= 
([q,\infty)\less[r,\infty))\cup([q,\infty)\less[s,\infty))$.  From the way in 
which we defined $K$ it follows that  
$\limsup\inorm{(f-f_n)\restriction _{[q,\infty)}}<2\epsilon $, which is what 
we want to prove. Since $|f(r)|\ge \epsilon $ for all $r\in H$, the set $H$ is 
relatively compact in $\Upsilon $ and so  has a supremum $s$ in $\Upsilon $, 
with $|f(s)|\ge \epsilon $ by continuity of $f$.  We shall next show that $s$ 
is in $H$.

To show that $\Phi (f_n;s)\to\Phi 
(f;s)$ is not hard, since by continuity of $\Phi (f;\cdot)$ we can 
find, given $\eta >0$, some $r\in H$ with $\Phi (f;r)<\Phi (f;s)+\eta $.
Since each $\Phi (f_n;\cdot)$ is a decreasing function, we therefore have
$$
\limsup\Phi (f_n;s) \le  \lim\Phi (f_n;r)  =  \Phi(f;r) < \Phi (f;s)+\eta .
$$

As in the proof of Lemma~6.4, we consider the finite set $M=\max\{u\in 
[s,\infty): |f(u)|\ge \epsilon \}$ and choose $m\in 
\nats$ such that  $|f(t)|<\epsilon -2^{-m}$ whenever
$t\in [s,\infty)\less\bigcup_{u\in M} [s,u]$. Given $\eta >0$ we may choose 
$r\in H$ such that $\Phi (f;r)<\Phi (f;s)+\eta $ and
$\inorm{f\restriction_{[r,s]}}<|f(s)|+\eta2^{-m} $. We may also suppose that
$\inorm{f\restriction_{[r,\infty)\less[r,s,\infty)}}<\epsilon $. Since 
$r\in H$ we have $\Phi (f_{n};r) \to \Psi (f;r)$ and so we may apply 
Lemma~6.4, obtaining $s'\in[r,\infty)$ such 
that 
$$
\align
|f(s')|+2^{-m}\Xi (f;r,s')&=\lim_{j\to\infty}\bigl(|f_{n}(s')|+
2^{-m}\Xi (f_{n};r,s')\bigr)\\
       &=\lim_{j\to\infty}\sup_{t\in [r,\infty)}\bigl(
|f_{n}(s)|+2^{-m}\Xi (f_{n};r,s)\bigr),
\endalign
$$
and such that, moreover, $\Phi (f_{n};s')\to \Phi (f;s')$. By our choice 
of $m$ it must be that $|f(s')|\ge \epsilon $. As in our proof that $s'$ is 
non-empty, we show that $s'$ is in $H$, which implies that $s'\trle s$. By the 
way in which we chose $r$, we have $|f(s')|\le |f(s)|+\eta $, so that
$$
\align
\limsup_{j\to \infty}\Xi (f_{n_{k_j}};r,s) &\le
\lim_{j\to\infty}\sup_{t\in [r,\infty)}\bigl(
2^m|f_{n}(s)|+\Xi (f_{n};r,s)\bigr)-2^m|f(s)|\\
  &= 2^m |f(s')|+2^{-m}\Xi (f;r,s')-2^m|f(s)|\\
  &\le \Xi (f;r,s') +\eta \\
  &\le \Xi (f;r,s) +\eta \\
  &\le\liminf_{j\to \infty}\Xi (f_{n};r,s)+\eta  .
\endalign
$$
Thus $\Xi (f_{n};r,s)\to \Xi (f;r,s)$, which, as we have seen, 
implies that $f_n$ converges to $f$ uniformly on $[r,s]$ and that 
$\inorm{f_{n}\restriction_{[r_1,\infty)\less[r_1,s,\infty)}}$ tends to 
$\inorm{f\restriction_{[r_1,\infty)\less[r_1,s,\infty)}}$, a quantity known 
to be smaller than $\epsilon $.  

We have now found a maximal element $s$ of $H$. To finish the proof we 
have to show that $\limsup\inorm{(f-f_n)\restriction _{[s,\infty)}}<2\epsilon 
$.

\proclaim{Claim 2}
Assume that one of the following holds:
\roster
\item"{\rm(a)}" $f$ is not monotone on $[s,\infty]$;
\item"{\rm(b)}" there is some $t\in [s,\infty)$ with $f(t)\ne \sum_{u\in F_t} 
f(u)$;
\item"{\rm(c)}" there is some $t\in \Upsilon ^+\cap (s,\infty)$ such that 
$\rho (t)\ne \rho (t^-)$ and $f(t)\ne0$.
\endroster
Then there exists $w\in (r,\infty)$ such that:
\roster
\item"{\rm (1)}" \ $f_{n}\to f$ uniformly on $[s,w]$;
\item"{\rm (2)}" \ $\Phi (f_{n}.\indic_{\Upsilon \less(0,w,\infty)};s)
     \to \Phi (f.\indic_{\Upsilon \less(0,w,\infty)};s)$;
\item"{\rm (3)}" \ $\Phi (f_{n};w)\to \Phi (f;w)$;
\endroster 
\endproclaim
\demo{Proof}
This is very similar to the proof of Lemma~6.4 except that in the cases (a), 
(b) and (c), we use the functions $\Psi $, $\Theta $ and $\Omega $, 
respectively, instead of $\Sigma $. Consider, for example, case (a), where $f$ 
is not monotone on $[s,\infty]$.  Since 
$|f(s)|=\inorm{f\restriction_{[s,\infty)}}$, we see that
$$
\sup_{s\trle t\trle u} |f(s)-f(t)+f(u)|>\inorm{f\restriction_{[s,\infty)}}.
$$
For suitably chosen $m\in\nats$ therefore, the supremum of $|f(s)-
f(t)+f(u)|+2^{-m}\Xi (f;s,t)$ is strictly greater than $\inorm
{f\restriction_{[s,\infty)}}+2^{-m}$. We apply Lemma~1.2 with $I=[s,\infty)$,
$\phi _t(f)=\sup_{u\in[t,\infty)} |f(s)-
f(t)+f(u)|+2^{-m}\Xi (f;s,t)$ and $\psi _t(f) = \Phi (f;t)$, getting a 
sequence $(t_n)$ which converges in $\Upsilon \cap\{\infty\}$. and which 
satisfies, among other things, $\phi _{t_n}(f)\to \sup_t\phi_t(f)$.  The 
limit cannot be $\infty$ and cannot be $s$ since in these cases we should have 
$\limsup\phi_{t_n}(f)\le \inorm{f\restriction_{[s,\infty)}}+2^{-m}$, which is 
a contradiction, since $(t_n)$ satisfies, among other things, $\phi 
_{t_n}(f)\to \sup_t\phi_t(f)$. We take $w=t$ and can now establish (1), (2) and 
(3) as in the proof of Lemma~6.4.  
\qed
\enddemo

Suppose now that we are in the situation where the above Claim is applicable 
and $w\in (s,\infty)$ has properties (1),(2),(3). If 
$\inorm{f\restriction_{[w,\infty)}} \ge \epsilon $ then $m(f.\indic_{\Upsilon 
\less(0,w,\infty)},s,\epsilon)$ is smaller than $m(f,q,\epsilon )$, so that 
$\limsup\inorm{(f-f_n)\restriction_{[s,\infty)\less[s,w,\infty)}} <2\epsilon $ 
by (2) and our inductive hypothesis. Taken together with (1), this tells us 
that $w$ is in $K$, so that $[w,\infty)\cap H\ne \emptyset$, contradicting the 
maximality of $s$ in $H$. 

Suppose now that 
$\inorm{f\restriction_{[s,\infty)\less[s,w,\infty)}}\ge\epsilon $.  Our 
inductive hypothesis, together with (3), gives $\limsup \inorm{(f-
f_n)\restriction_{[w,\infty)}}<2\epsilon $.  It is not hard to see that if we  
apply Lemma~6.4 to the functions $f_n.\indic_{\Upsilon \less(0,w,\infty)}$ and 
$f.\indic_{\Upsilon \less(0,w,\infty)}$ then we obtain an element of 
$H\cap(s,\infty)$, again contradicting maximality of $s$.  

The only remaining possibility is that both
$\inorm{f\restriction_{[w,\infty)}} $ and
$\inorm{f\restriction_{[s,\infty)\less[s,w,\infty)}}$ are smaller than 
$\epsilon $.  In this case, the initial ``$m(f,q,\epsilon)=0$'' case of the 
proposition gives us $\limsup \inorm{f\restriction_{[w,\infty)}}<2\epsilon $ 
and $\limsup\inorm{f\restriction_{[s,\infty)\less[s,w,\infty)}}<2\epsilon $. 
Taken together with (1), these give us what we want to prove. 

Finally we have to deal with the case where none of the possibilities
(a),(b),(c) occurs. 

\proclaim{Claim 3}
If none of {\rm(a),(b),(c)} holds then:
\roster
\item $f$ is a monotone function on $[s,\infty]$;
\item for every $t\in [s,\infty)$ we have $f(t)=\sum_{u\in F_t} f(u)$; 
\item $\supp f\cap [s,\infty)\subseteq \{t\in [s,\infty):\rho (t)=\rho (s)\}$;
\item for each $l\ge 1$ and each $t\in [s,\infty)$, $\Alpha_l(f,t)
=l^{-1}|f(t)|$;
\item for any $\eta >0$ there exists a finite antichain $A\subseteq 
[s,\infty)$ such that $0<|f(t)|<\eta $ for all $t\in A$ and such that 
$\sum_{t\in A}|f(t)|>|f(s)|-\eta $.
\endroster
\endproclaim
\demo{Proof}
Assertions (1) and (2) follow immediately from the failure of (a) and (b). 
If (3) were not true, there would exist $t$ minimal subject to $f(t)\ne0$ and 
$\rho (t)>\rho (s)$; by monotonicity, $t$ would be minimal subject to $\rho 
(t)>\rho (s)$ and so $t$ would be in $\Upsilon ^+$ by the assumed continuity 
of $\rho $. Thus (c) would hold, contrary to supposition.  

Without loss of generality, we may suppose that $f(s)>0$. The monotonicity now
implies that $f$ is everywhere non-negative on $[s,\infty)$, since by 
convention $f(\infty)$ has been taken to be $0$. To prove (4) we proceed by 
induction on $l$, the case $l=1$ being true by monotonicity. Now let $t$ be 
in $[s,\infty)$ and let $A\subseteq 
[t,\infty)$ be an antichain of size $l$ with $f(w)>0$ for all $w\in A$. 
Define $u$ to be the greatest element such that $u\trle w$ for all $w\in 
A$. By (3), $\rho (w)=\rho (s)$ for all $w\in A$, and so $A=\bigcup_{v\in 
F_u} A\cap [v,\infty)$.  By choice of $u$ we have $\#(A\cap [v,\infty))<l$ for  
all $v$, whence $\sum_{w\in a\cap [v,\infty)}f(w)\le f(v)$ by inductive 
hypothesis.  Thus
$$
\sum_{w\in A} f(w)= \sum_{v\in F_u}\sum_{v\in A\cap[u,\infty)}f(w)
\le \sum_{v\in F_u}f(v) = f(u)\le f(t)
$$
by (2) and monotonicity.

To prove (5), we consider the infinite antichain $B=\min\{v\in 
[s,\infty):0<f(v)<\eta \}$; provided we can prove that 
$f(t)=\sum_{v\in B}f(v)$, it is clear that a suitable finite subset of $B$ 
will do for $A$. Define $g(t)=\sum_{v\in B\cap(t,\infty)}f(v)$ and consider 
$g$ as a function on the set $T= \{t\in [s,\infty):f(t)\ge \eta \}$. Notice 
that $g$ is continuous and that for all $t\in T$ we have $g(t) = \sum_{u\in 
t^+\less T} f(t) + \sum_{u\in t^+\cap T} g(t)$. Suppose if possible that the 
subset $U=\{u\in T:f(u)\ne g(u)\}$ is not empty. The subset $U$ can have no 
maximal elements since if $u$ were such an element we should have 
$$
g(u)= \sum_{v\in u^+\less T} f(v) + \sum_{v\in u^+\cap T} g(v)
    = \sum_{v\in u^+} f(v) = f(u),
$$
by (2), (3) and maximality of $u$ in $U$. Thus there exists some $w\in T\less 
U$ which is a limit point of $U$.  For some $t\in  U\cap [s,w)$ the set $T\cap 
[t,w]$ is totally ordered; equivalently, for each $u\in [t,w)$ there is 
a unique $u'\in u+$ with $f(u')\ge \eta $.  We have 
$$
g(u)-g(u') = \sum_{v\in u^+\less T} f(v) = f(u)-f(u'),
$$
which, with the continuity of $g$ and $f$, implies that $g(t)-g(w) = f(t)-
f(w)$. But this is a contradiction, since we had $f(w)=g(w)$ and $f(t)\ne 
g(t)$.
\qed\enddemo

We are at last near the end of the proof of Proposition~6.5.  The reader will 
recall that we are trying to show that $\limsup\inorm{(f-
f_n)\restriction_{[s,\infty)}}<2\epsilon $.  We take $\eta =\quarter\epsilon $ 
and get a set $A$ as in (5); let $\#A=l$. Let $k$ be the integer part of 
$f(s)/\eta $ and   for each $w\in A$ and each natural 
number $j\le k+1 $, define 
$$
v(w,j)=\max\{v\in [s,w]:f(v)\ge f(s)-j\eta \}.
$$ 
Notice that $v(w,k+1)=w$.
Let $u(w,0)=s$ and for $1\le j\le k+1<f(s)/\eta $ let $u(w,j)$ be the unique 
element of $(s,w]\cap v(w,j-1)^+.$ 

Since $\Phi (f_n;s)\to \Phi (f;s)$ we have $\Delta ^+(f_n;s)\to \Delta 
^+(f;s)=0$ and $A_l(f_n;s)\to A_l(f;s)$ which equals $l^{-1}f(s)$. Remembering 
that $f_n$ converges pointwise to $f$ we have, for all sufficiently large $n$,
$$
\Delta ^+(f_n;s)<\eta ,\text{ and } A_{l+1}(f_n;s)< f(s)+(l+1)^{-1}\eta,
$$
as well as $|f_n(x)-f(x)|<\eta $ for all $x$ in the finite set 
$\bigcup_{j,w}\{u(w,j),v(w,j)\}$ and $|f_n(w)-f(w)|<\eta/l $ for all $w\in A$. 
We show that $|f(t)-f_n(t)|<2\epsilon $ for any such $n$ and any $t\in 
[s,\infty)$.  There are three cases to consider. 

First suppose that $t$ is in $[w,\infty)$ for some $w\in A$.  Since $\Delta 
^+(f_n)<\eta $ there is a decreasing function $g$ such that $\inorm{(g-
f_n)\restriction_{[s,\infty)}}<\eta $.  Hence,  $-\eta <f_n(t)<f_n(w)+2\eta 
<f(w)+3\eta <4\eta $. Since also $0\le f(t)<\eta $ we have $|f_n(t)-
f(t)|<4\eta $.

Next suppose that $t$ is in $[s,w]$ for some $w\in A$.  Then, for a suitably
chosen $j$ we have $u(w,j)\trle t\trle v(w,j)$. So, again using the fact that 
$f_n$ is within $\eta $ of  a decreasing function, we have
$f_n(t)< f_n(u(w,j)+2\eta <f(u(w,j))+3\eta \le f(v(w,j))+4\eta 
<f_n(v(w,j))+5\eta <f_n(t)+7\eta $, whence again $|f_n(t)-f(t)|<2\epsilon .$

Finally, suppose that $t$ is incomparable with all $w\in A$.  Since 
$A\cup\{t\}$ is an antichain of size $l+1$, we have 
$f_n(t)+\sum_{w\in A}f_n(w)\le (l+1)A_{l+1}(f_n;s)<f(s)+\eta $. But
$\sum_{w\in A}f_n(w)>\sum_{w\in A} f(w) - \eta >f(s)-2\eta $, whence
$f_n(t)<3\eta $.  Since also $-\eta <f_n(t)$ and $0\le f(t)<\eta $, we have 
$|f_n(t)-f(t)|<3\eta <2\epsilon $ in this case too.

\proclaim{Corollary 6.6}
There is an Asplund space $X$ which admits a Kadec renorming but no strictly 
convex renorming.
\endproclaim
\demo{Proof}
If $\Upsilon $ is a full dyadic tree of height $\omega _1$ then any increasing 
function $\rho :\Upsilon \to \reals$ is constant on some set $[u,\infty)$.  
Hence, by 3.4, $\C_0(\Upsilon )$ is not strictly convexifiable. On the other 
hand this space admits a Kadec renorming, since all points of a finitely 
branching tree are good, even for a constant function $\rho $.
\qed\enddemo

As promised earlier, we finish this section with a result about 
$\sigma $-fragmentability.

\proclaim{Proposition 6.7}
For a tree $\Upsilon $ the space $\C_0(\Upsilon )$ admits a Kadec renorming if 
and only if it is $\sigma $-fragmentable.
\endproclaim
\demo{Proof}
As we have already remarked, one implication is already known.  Let us assume 
therefore that $\C_0(\Upsilon )$ is $\sigma $-fragmentable.  By the results of 
[\Haydonscat] $\Upsilon$ is a countable union $\bigcup_{n\in \nats}\Delta _n$ 
of subsets $\Delta _n$ that are discrete in the reverse topology.  For $t\in 
\Upsilon $ we define $M_t=\{m\in \nats: [t,\infty)\cap \Delta _m=\emptyset\}$ 
and set $\rho (t)=\sum_{m\in M_t}2^{-m}$. Obviously $\rho $ is an increasing 
function. We shall show that all points $t$ are $\rho $-good.  Any $t$ is in 
$\Delta _n$ for some $n$ and, by reverse-discreteness, there is a reverse-open 
neighbourhood $U$ of $t$ such that $U\cap \Delta _n=\{t\}$.  We may take $U$ 
to be of the form $[t,\infty)\less\bigcup_{u\in F}[u,\infty)$ for a suitable 
finite subset $F$ of $t^+$.  Evidently, $\rho (v)\ge \rho (t)+2^{-n}$ for all 
$v\in t^+\less F$, which shows that $t$ is a good point for $\rho $. \qed 
\enddemo

\heading 7. Strictly convex dual norms \endheading

In this section we present a sufficient condition for the existence on 
$\C_0(\Upsilon )$ of an equivalent norm having strictly convex dual norm. 
The reader should be warned that this condition is probably quite 
far from being necessary, though we shall see in Section~10 that it is 
sufficiently general to be satisfied even in some cases where $\C_0(\Upsilon 
)$ has neither a Fr\'echet-differentiable renorming nor a strictly convex 
renorming. 

\proclaim{Theorem 7.1}
Suppose that on the tree $\Upsilon $ there is an increasing function $\rho 
:\Upsilon \to \reals$ which is constant on no strictly increasing sequence in 
$\Upsilon$. Then there is an equivalent norm on $\C_0(\Upsilon )$ with 
strictly convex dual norm. 
\endproclaim 
\demo{Proof}
The dual space of $\C_0(\Upsilon )$ may be identified with $\ell_1(\Upsilon)$ 
and the norm dual to the supremum norm is of course the usual $\ell_1$-norm
$\norm{\xi }_1=\sum_{t\in \Upsilon }|\xi (t)|$. We shall construct an 
equivalent norm on $\ell_1(\Upsilon )$ which is strictly convex and lower 
semicontinuous for the weak* topology $\sigma (\ell_1(\Upsilon ),\C_0(\Upsilon 
))$. Useful weak* lower semicontinuous functions to use as building blocks for 
our norm are functions of the form $\xi \mapsto \norm{\xi \restriction V}_1$ 
with $V$ an open subset of $\Upsilon $. In particular, we may take $V$ of the 
form $[s,\infty)$ or $\{s\}$ when $s$ is a successor element of $\Upsilon $.

Let us write $\Upsilon _0$ for the set of all $t\in \Upsilon^+$ such that $\rho 
(t)>\rho (t^-)$. Notice that without spoiling the assumed property of 
$\rho $ we may modify that function so that it takes rational values at all 
points of $\Upsilon _0$. 

For each rational $q$ we note that the wedges $[s,\infty)$, with 
$s\in \Upsilon_0$ and $\rho (s)=q$ are disjoint, so that the family 
$(\norm{\xi \restriction[s,\infty)}_1)_{s\in\Upsilon_0\cap\rho^{-1}(q)}$ is 
in $\ell_1(\Upsilon_0\cap\rho^{-1}(q))\subseteq 
c_0(\Upsilon^+\cap\rho^{-1}(q))$. If we define 
$$ 
\align 
\Phi(q;\xi) &= \norm{(\norm{\xi\restriction[s,\infty)}_1)_{s\in \Upsilon_0 \cap 
\rho ^{-1}(q)}}_{\text{Day}}, \\ 
&=\bigl[\sum_{m=1}^\infty 2^{-m}\sup_{F\subset\Upsilon_0\cap\rho^{-1}(q),\ 
\#F=m}\sum_{s\in F} \norm{\xi\restriction[s,\infty)}_1^2\bigr]^{\frac12},
\endalign 
$$ 
the function $\Phi (q;\cdot)$ is weak* lower semicontinuous by our earlier 
observation. So  also is the norm defined (using suitable positive constants 
$c(q)$) by
$$
\norm{\xi }^2 = \norm{\xi }_1^2 + \norm{(|\xi (s)|)_{s\in \Upsilon 
^+}}_{\text{Day}}^2 + \sum_{q\in\rats} c(q) \Phi (q;\xi )^2.
$$

To see that this norm is strictly convex notice first that if $\xi ,\eta \in 
\ell_1(\Upsilon )$ are such that $\norm{\xi }=\norm{\eta }=\half\norm{\xi 
+\eta }$ then $\xi (s)=\eta (s)$ for all $s\in \Upsilon ^+$ because of the 
second term in the definition of $\norm{\cdot}$. Next, because of the $\Phi $ 
terms, we must have 
$$
\norm{\xi \restriction[s,\infty)}_1
=\norm{\eta \restriction[s,\infty)}_1
=\half\norm{(\xi+\eta) \restriction[s,\infty)}_1
$$ 
for every $s\in \Upsilon _0$.

We now consider a limit element $t$ of $\Upsilon $. Since $\rho $ is constant 
on no strictly increasing sequence in $\Upsilon $, there exist elements $s_n$ 
of $\Upsilon_0$ such that $t=\lim s_n$. The sets $[s_n,\infty)$ form a 
decreasing sequence with intersection $[t,\infty)$ and so we can now conclude 
that
$$
\norm{\xi \restriction[t,\infty)}_1
=\norm{\eta \restriction[t,\infty)}_1
=\half\norm{(\xi+\eta) \restriction[t,\infty)}_1\,.
$$ 
Again using our assumption about $\rho $ we see that the minimal elements of 
$\{v\in (t,\infty): \rho (v)>\rho (t)\}$ are all in $\Upsilon _0$ and that
$\{u\in (t,\infty):\rho (u)=\rho (t)\}$ contains only points of $\Upsilon ^+$.
We have
$$
\norm{\xi \restriction[t,\infty)}_1 = |\xi (t)|+\sum_{u\in(t,\infty),\ \rho 
(u)=\rho (t)}|\xi (u)| + \sum_{u\in\min\{v\in(t,\infty): \rho (v)>\rho (t)\}}
\norm{\xi \restriction[v,\infty)}_1,
$$
with similar expressions for $\eta $ and $\xi +\eta $. From what we have 
already shown we can now deduce that
$$
|\xi (t)|=|\eta (t)|=\half|\xi (t)+\eta (t)|,
$$
which is to say $\xi (t)=\eta (t)$.
\qed
\enddemo

\heading 
8. Fr\'echet-differentiable renormings 
\endheading

In this section we shall establish a necessary and sufficient condition for 
the existence on $\C_0(\Upsilon )$ of an equivalent Fr\'echet-differentiable 
norm. Rather surprisingly, this condition turns out to be the same as the one 
obtained in Section~4 for LUR renormings. 

\proclaim{Theorem 8.1} 
For a tree $\Upsilon $ the following are equivalent:
\roster
\item there is a $\C^\infty$ renorming of $\C_0(\Upsilon )$;
\item there is a Fr\'echet-differentiable renorming of $\C_0(\Upsilon )$;
\item there is an increasing real-valued function $\rho $ on $\Upsilon $ which 
has no bad points and is constant on no ever-branching subset;
\item $\C_0(\Upsilon )$ admits a Talagrand operator.
\endroster
\endproclaim

It is only the implications (2)$\implies$(3) and (3)$\implies$(4) which need 
proof and we shall deal with the second of these first. We start by making 
an observation about sets with no ever-branching subset.

\proclaim{Lemma 8.2}
Let $\Upsilon $ be a tree and let $U\subseteq\Upsilon $ have no
ever-branching subset. We may define an ordinal-valued function $i_U$ on $U$ 
 which has the following properties:
\roster
\item $i_U$ is decreasing for the tree-order;
\item for any totally-ordered subset $V$ of $U$, $i_U$ takes only finitely 
many values on $V$;
\item for each $t\in U$ there is at most one $u$ in $t^+\cap U$ with 
$i_U(u)=i_U(t)$.
\endroster
\endproclaim
\demo{Proof}
We define $i_U$ in a standard way, starting from a derivation. For a non-empty 
subset $W$ of $U$ we define 
$$
W'= W\less\{w\in W: W\cap [w,\infty) \text{ is totally ordered} \}.
$$
Since $U$ has no ever-branching subset, $W'$ is a proper subset of $W$ 
whenever $W$ is non-empty. We then define recursively
$$
\align
U^{(0)} &= U, \\
U^{(\beta )} &= \bigcap_{\alpha <\beta } (U^{(\alpha )})',
\endalign
$$
and set $i_U(t)= \beta $ when $t\in U^{(\beta )}\less U^{(\beta +1)}$.

It is clear that $i_U$ is decreasing for the tree-order. If $V$ is a totally 
ordered subset of $U$ then $V$ is well-ordered and a decreasing ordinal-valued 
function on a well-ordered set can take only finitely many values. Finally, if 
$t$ is in $U$ and $i_U(t)=\beta $ then $U^{(\beta )}\cap [t,\infty) $ is 
totally ordered, so that $U^{(\beta )}$ can contain at most one element of 
$t^+$. \qed\enddemo 

Let us now suppose that $\Upsilon $ is a tree equipped with an increasing 
real-valued function $\rho $ which has no bad points and which is constant on 
no ever-branching subset. As usual, we may suppose that 
$\rho $ takes values in $(0,1)$ and is continuous for the locally compact 
topology of $\Upsilon $. We shall say that a pair $(s,u)$ in $\Upsilon 
\times\Upsilon $ is a {\sl special pair } if either $s=u$ or $\rho (s)=\rho 
(u)$ and there exists $t\in s^+$ with $t\trle u$ and $i_{\rho ^{-1}(\rho 
(s))}(t)<i_{\rho ^{-1}(\rho (s))}(s)$. It follows from assertion (2) of Lemma 
8.2 that for any $u\in \Upsilon $ there are only finitely many $s$ such that 
$(s,u)$ is special. We define $T:\C_0(\Upsilon )\to c_0(\Upsilon \times 
\Upsilon )$ by setting 
$$
(Tf)(s,u) = \cases
\displaystyle\frac{\delta _u}{1+\#F_u}\biggl[f(u)-\sum_{v\in F_u} f(v)\biggr] 
            \quad\text{ if $(s,u)$ is special},\\
            0 \qquad\qquad\qquad\qquad\qquad\qquad \text{otherwise},
            \endcases
$$
where, as previously, we write $F_u$ for the finite set $ \{v\in u^+:\rho 
(v)=\rho (u)\}$.

Clearly $T$ is a bounded linear operator from $\C_0(\Upsilon )$ into 
$\ell_\infty(\Upsilon \times\Upsilon )$ with $\norm{T}\le 1$. To show that $T$ 
takes values in $c_0(\Upsilon \times\Upsilon )$ we note that when $(Tf)(s,u)$ 
is not zero it equals $(Sf)(u)$, where $S:\C_0(\Upsilon )\to c_0(\Upsilon )$ 
is the operator defined  in Section 4. For any $f\in \C_0(\Upsilon )$ and any 
$\epsilon >0$, Lemma~4.3 tells us that there are only finitely many $u\in 
\Upsilon $ with $|(Sf)(u)|>\epsilon $. For each of these $u$ there are only 
finitely many $s$ with $(s,u)$ special and so $|(Tf)(s,u)|$ exceeds $\epsilon 
$ for only finitely many pairs $(s,u)$. We have now shown that $T$ takes 
values in $c_0(\Upsilon \times \Upsilon )$. 

To show that $T$ is a Talagrand operator we let $f$ be a non-zero element of 
$C_0(\Upsilon )$ and let $s$ be an element of $\Upsilon $ maximal in the tree 
order subject to the condition $|f(s)|=\inorm{f}$. If $f(s)\ne \sum_{t\in F_s} 
f(t)$ then we have $(Tf)(s,s)\ne 0$. Suppose then that $f(s)=
\sum_{t\in F_s} f(t) $. By assertion (3) of Lemma 8.2, there is at most one 
element $t_0$ of $F_s$ with $i_{\rho ^{-1}(\rho (s))}(t_0)=
i_{\rho ^{-1}(\rho (s))}(s)$. By maximality, we certainly do not have 
$f(s)=f(t_0)$ for such a $t_0$, and so there must exist some $t\in 
F_s$ with $i_{\rho ^{-1}(\rho (s))}(t)<i_{\rho ^{-1}(\rho (s))}(s)$ and 
$f(t)\ne 0$. Now by Lemma~4.4 there exists some $u\trge t$ with $\rho 
(u)=\rho (t)=\rho (s)$ and $(Sf)(u)\ne 0$. For this $u$, the pair $(s,u)$ is 
special and so $(Tf)(s,u)=(Sf)(u)\ne 0$. This completes the proof that that 
(3) implies (4) in the theorem.

\smallskip

We now move on to the implication (2)$\implies$(3). Let us suppose that 
$\norm{\cdot}$ is a Fr\'echet-smooth norm on $\C_0(\Upsilon )$,
equivalent to $\inorm{\cdot}$. We may assume that $\norm{\cdot}\le 
\inorm{\cdot} \le M\norm{\cdot}$, where $M$ is a positive constant. 
We shall construct an increasing sequence of increasing functions $\rho 
_n:\Upsilon \to \reals$ and an increasing sequence of subsets $\Delta _n$ of 
$\Upsilon $. To specify precisely the properties possessed by our system,
it will be convenient to introduce some further notation.  We 
write $\Beta_n$ for the set of all $t\in \Upsilon $ which are bad points or 
fan points for $\rho _n$. For $t\in \Beta _n$ we set 
$$
\epsilon^n _t= \left\{ \matrix \format \l & \qquad \l \\
      \rho _n(t) & \text {if $t$ is a minimal element of $\Beta_n$}\\
      \rho_n(t)-\sup\{\rho _n(s):s\in (0,t)\cap \Beta_n\} &\text{otherwise}
                       \endmatrix\right.
$$ 
and let $\Gamma _n$ be the set of all $t$ in $\Beta_n$ with $\epsilon ^n_t>0$.
We notice that 
$$
\rho _n(t) = \sup\{\rho _n(s):s\in (0,t)\cap \Gamma _n\}
$$
whenever $t\in \Beta_n\less\Gamma_n$. 

For $t\in \Upsilon $ and $n\in \nats$, we define $C^n_t$ be the set of 
functions $f\in \C_0(\Upsilon )$ such that 
\roster
\item $\supp f=(0,t]$;
\item $f\restriction (0,t]$ is increasing;
\item whenever $r\in (0,t)$ and $s$ is the unique member of $t^+\cap(0,t]$, we 
have $f(r)=f(s)$ unless $r\in \Delta _n$.
\endroster
The content of condition (3) is that, regarded as increasing functions on 
$(0,t]$, the functions in $C^n_t$ have jumps only at points of $\Delta _n$.  

Finally, for all $t$, we define
$$
\delta ^n_t= \inf_{u\in t^+}\rho _n(u)-\rho _n(t).
$$

The crucial properties of the system to be constructed are:
\roster
\item"(a)" $\rho _n$ takes values in the interval $(0,2-2^{-n})$;
\item"(b)" $\delta _n(t)>0$ for all $t\in \Delta _n$;
\item"(c)" for $n\ge 1$, $\bigcup_{m<n}\Gamma _m\subseteq \Delta_n$;
\item"(d)" if $u\in \Beta_n$ and $g\in C^n_u$ then $u$ is a bad point or a fan 
point for the function $\mu (g,\cdot)$;
\item"(e)" if $u\in\Beta_n\less \Gamma _n$ and $g\in C^n_u$ then $\mu (g,u) = 
\sup\{\mu(g.\indic_{(0,t]},t):t\in \Gamma _n\cap (0,u)\}$.
\endroster
We now show how to perform the construction so that (a) to (e) hold.

We start by setting $\rho _0(t)=\mu(t) $ and $\Delta _0=\emptyset$. Thus (a), 
(b) and (c) certainly hold for $n=0$. Applying our definitions, $B_0$ 
is the set of bad points and fan points for the $\mu $ function.  Now, since 
$\Delta _0=\emptyset$, each set $C^n_t$ contains only scalar multiples of the 
indicator function $\indic_{(0,t]}$.  For such a function $g=\lambda 
\indic_{(0,t]}$, we have $\mu (g,u)=|\lambda| \mu (u)=|\lambda|\rho _0(u)$ for 
all $u\in [t,\infty)$, and so (d) holds as well. The final assertion (e) 
follows from the observation we made just after the definition of $\Gamma _n$.

Now let us suppose that $\rho _0, \dots,\rho _n$ and $\Delta_0,\dots,\Delta 
_n$ have been constructed and satisfy (a) to (e).  It follows from (b) that 
$\Delta _n\cap (0,t)$ is countable for each $t\in \Upsilon $ (indeed, 
$\sum_{s\in \Delta _n\cap(0,t)}\delta ^n_s\le \rho (t)$).  So each set $C^n_t$ 
is norm-separable; we fix, for each $t$, a norm-dense sequence 
$(f^{n,t}_k)_{k\in\nats}$ in this set and also fix an enumeration 
$(q_l)_{l\in\nats}$ of the non-negative rationals. We define 
$$ 
\sigma^n _t(u) =\left\{\matrix \format \l &\qquad\l \\ 
\sum_{k,l\in\nats} 2^{-k-l-2}(\inorm{f^{n,t}_k}+q_l)^{-1}\mu (f^{n,t}_k,q_l,u)& 
                                          \text{if $t\trl u$}\\
         0&\text{otherwise,}
              \endmatrix\right.
$$
where 
$$
\align
\mu (f,\delta ,u)&= \mu(f+(f(t)+\delta )\indic_{(t,u]},u)\\
                 &= \inf\{\norm{f+(f(t)+\delta )\indic_{(t,u]}+g}: 
g\in\C_0(\Upsilon ) \text{ and } \supp g\subseteq (u,\infty)\}
\endalign
$$
is another variation on a familiar theme.
We notice that $\sigma _t(u)\le 1$ and that $\sigma _t(u)\ge M^{-1}$ when 
$u\in (t,\infty)$.

We now define $\Delta _{n+1}=\Delta _n\cup\Gamma_n$ and set
$$
\rho _{n+1} = \rho _n + 2^{-n-3}\sum_{t\in \Gamma _n} \epsilon^n _t \sigma^n 
_t + 2^{-n-3}\sum_{t\in\Delta _n} \delta ^n_t\sigma^n_t.
$$
Of course, we need to show that the series defining $\rho _{n+1}$ converges and 
that $\rho _{n+1}$ takes values in $(0,2-2^{-n-1})$ in order that (a) should 
hold at stage $n+1$. This is indeed so since for any $u\in\Upsilon $ we have 
$$
\align
\rho_{n+1}(u) &\le \rho_n(u)+2^{-n-3} \sum_{t\in \Gamma_n\cap(0,u)}\epsilon^n_t
 + 2^{n-3} \sum_{t\in \Delta_n\cap(0,u)}\delta^n_t\\
& \le (1+2.2^{-n-2})\rho_n(u)\le (1+2^{-n-2})(2-2^{-n}) < 2-2^{-n-1}.
\endalign
$$
To see that (b) holds we note that if $t\in \Gamma _n$ and $u\in t^+$ then 
$\sigma ^n_t(t)=0$ and $\sigma ^n_t(u)\ge 1/M$ so that $\delta ^{n+1}(t)\ge
2^{n-3}\epsilon ^n_t/M>0$. Assertion (c) holds because of the way we defined 
$\Delta _{n+1}$.

Now let $u$ be a bad point for $\rho _{n+1}$.  It is clear from the definition 
that if $t\in \Delta _n\cap(0,u)$ and $k,l$ are natural numbers then $u$ is 
bad for the function $\mu (f^{n,t}_k,q^l,\cdot)$.  By uniform 
approximation, we see that the same is true for all functions $\mu (f,\delta 
,\cdot)$ with $f\in C^n_t$ and $\delta \in\reals^+$. Looking at this another 
way, we see that $u$ is bad for each function $\mu (g,\cdot)$ where $g$ is a 
function in $C^n_u$ of the form $f+(f(t)+\delta )\indic_{(t,u]}$ with $t\in 
\Delta _n\cap(0,u)$, $f\in C^n_t$ and $\delta \ge 0$.  The set of functions of 
this form being norm-dense in $C^n_u$, we see that $u$ is bad for all of the 
functions $\mu (g,\cdot)$ with $g\in C^n_u$.  A similar argument works for the 
case where $u$ is a fan point for $\rho _{n+1}$, since if $\rho_{n+1}(v)=\rho 
_{n+1}(u)$ for all $v$ in some ever-branching subset $T$ of $[u,\infty)$ 
then all the functions $\mu (f^{n,t}_k,q^l,\cdot)$ are constant on the same 
set $T$.       In this way we establish (d).

We have already noted that 
$$
\rho_{n+1} (u)=\sup \{\rho _{n+1}(t): t\in \Gamma _{n+1}\cap (0,u)\}
$$
whenever $u\in \Beta_{n+1}\less\Gamma _{n+1}$.  It follows from the 
definition of $\rho _{n+1}$ that 
$$
\mu (f^{n,t}_k,q_l,u)=\sup \mu (f^{n,t}_k,q_l,t): t\in \Gamma _{n+1}\cap (0,u)\}
$$
for all $t\in \Delta _n$ and all $k,l$.  A uniform approximation argument now 
leads us to (e).

Having completed the construction of the sequence $(\rho _n)$ we define 
$\rho (t)$ to be $\lim_{n\to\infty}\rho _n(t)$.  It is this function that we 
shall show has no bad points and no fan points.  Suppose then that $u$ is a 
bad point or a fan point for $\rho $.  By construction, $u$ is in $\Beta_n$ 
for all $n$, and since $\Beta_{n+1}\cap \Gamma_{n}\subseteq 
\Beta_{n+1}\cap\Delta _{n+1}=\emptyset$ by (b) and (c), $u$ is in none of the 
sets $\Gamma _n$.  

We shall now construct an increasing sequence $(t_n)$ in $(0,u)$, with $t_n\in 
\Gamma _n$, and a summable sequence $(\delta _n)$ of non-negative real 
numbers, and shall define $f_n\in C^n_{t_n}$, $g_n\in C^n_u$ by
$$
\align
f_n & = \indic_{(0,t_n]} + \sum_{m=0}^{n-1} \delta _m\indic_{(t_m,t_n]},\\
g_n & = \indic_{(0,u]}   + \sum_{m=0}^{n-1} \delta _m\indic_{(t_m,u]}.
\endalign
$$
Since $u$ and $t_n$ are bad or fan for $\mu (f_n,\cdot)$, we may consider, as 
in Proposition~3.5, $\mu $-attaining modifications $\hat f_n$ and $\hat g_n$ of 
$f_n$.  The construction will be carried out in such a way that the following 
hold: 
\roster
\item"(i)" for all $n$, $\mu (f_n,t_n)\ge \mu (f_n,u)-2^{-n}$;
\item"(ii)" for each $n$, $\<\hat 
f_{n+1}\dual,\indic_{(t_n,t_{n+1}]}\>\ge1/4M$.
\endroster
Before giving the details, let us indicate why we shall end up with a 
contradiction.  

Define the function $g$ by  
$$
g=\indic_{(0,u]}+\sum_{m=0}^\infty \delta _m\indic_{(t_m,u]} .
$$
Thus $g$ is the norm limit as $n\to \infty$ of $g_n$. So by (d) and uniform 
approximation, $u$ is a bad point or a fan point for $\mu (g,\cdot)$.  
By (i), we have 
$$
\align
\norm{\hat g} &= \mu (g,u)\\
              &= \lim_{n\to\infty} \mu (g_n,u) \\
              &= \lim_{n\to\infty} \mu (f_n,t_n) \\
              &= \lim_{n\to\infty} \norm{\hat f_n}.
\endalign
$$
Also $\supp(g-f_n)\subseteq (t_n,\infty)$ and $\supp \hat f_n\dual \cap 
(t_n,\infty)=\emptyset$ by Lemma~3.5.  Thus
$$
\<\hat f_n\dual,\hat g\> = \<\hat f_n\dual, f_n\> = \norm{\hat f_n}
$$
which tends to $\norm{\hat g}$ as $n\to \infty$.  Hence $\norm{\hat g\dual -
\hat f_n\dual}\dual \to 0$ by \v Smulyan's Criterion.  But of course this 
contradicts (ii).

To finish the proof we have to show how to construct $t_n$ and $\delta _n$.
We start by choosing an arbitrary $t_0\in \Gamma _0\cap (0,u)$. (There are such
points since $u$ is in $\Beta_0\less \Gamma_0$.) If
$t_0\trl t_1\trl\cdots\trl t_n$ and $\delta _0,\dots,\delta _{n-1}$ have been 
defined already (so that $f_n\in C_{t_n}$ is determined as above), we choose 
$\delta_n$ to be the greatest $\delta \ge 0$ such that 
$$
\mu (f_n,\delta,u)\le \mu (f_n,0,u)+\delta/2M .
$$
To see that this definition makes sense, we note that the set of $\delta $ 
satisfying this inequality is non-empty (since it certainly contains $0$) 
and bounded above, since we have $\mu (f_n,\delta,u)\ge \delta/M$ for any 
$\delta$.  Now, by construction, $u$ is a bad or fan point for each of the 
functions $\mu (f_n,\delta ,\cdot)$ and so, as in Proposition~3.5, it is the 
case that for all $\delta \ge 0$
$$ 
\mu (f_n,\delta ,u)=\norm{f_n+(f_n(t_n)+\delta )(\indic_{(t_n,u]}+\phi _u)}, 
$$                  
where $\phi _u$ is either 0 or a fan function at $u$.  Thus $\delta _n$ is 
maximal subject to 
$$
\norm{f_n+(f_n(t_n)+\delta_n )(\indic_{(t_n,u]}+\phi _u)} \le
\norm{f_n+f_n(t_n)(\indic_{(t_n,u]}+\phi _u)} + \delta /2M. 
$$                  
We recall that the function 
$f_n+(f_n(t_n)+\delta_n)(\indic_{(t_n,u]}+\phi _u)$ is what we have decided 
to call $\hat g_{n+1}$.  Because of the maximality in the definition of $\delta 
_n$ it must be that $\<\hat g_{n+1}\dual, \indic_{(t_n,u]}+\phi _u\>\ge 1/2M$.  

The function $g_n$ is in $C^{n+1}_u$ and $u$ is in $\Beta_{n+1}\less \Gamma 
_{n+1}$.  So by (e) there is a sequence $(s_k)$ in $\Gamma_{n+1}\cap(0,u)$ 
with 
$$
\mu (g_n.\indic_{(0,s_k]},s_k)\to \mu(g_n,u)
$$
as $n\to \infty$.  We shall choose $t_{n+1}$ to be one of the $s_k$; 
evidently, any sufficiently large $k$ will give us (i). 
If we write $h_k=g_n.\indic_{(0,s_k]}$ then $s_k$, 
being an element of $\Gamma _{n+1}$, is a bad point or a fan point for $\mu 
(h_k,\cdot)$ so that we can define $\hat h_k$ as usual.  By an argument we used 
earlier, we have
$$
\<\hat h_k\dual,\hat g_n\> = \<\hat h_k\dual, h_k\> = \norm{\hat h_k} 
=\mu(h_k,s_k) \to \mu (g_n,u) = \hat g_n,
$$
so that $\norm{\hat h_k\dual - \hat g_n\dual}\dual \to 0$ by \v Smulyan's 
Criterion.  In particular, for sufficiently large $k$, we have 
 $\<\hat h_k\dual, \indic_{(t_n,u]}+\phi _u\>\ge 1/4M$. Since $\supp \hat 
h_k \cap (s_k,\infty)=\emptyset$ we have, in fact, 
$$
\<\hat h_k,\indic_{(t_n,s_k]}\> \ge 1/4M,
$$

All that remains is to verify that the series $\sum_{n=0}^\infty \delta _n$ 
converges.  To do this, we notice first that, by the way in which $\delta _n$ 
was chosen, we always have
$$
\mu (f_n,\delta_n,u)= \mu (f_n,0,u)+\delta_n/2M,
$$
or, equivalently,
$$
\mu (g_{n+1},u) = \mu (g_n,u)+\delta_n/2M.
$$
Thus, for any $N\ge 1$,
$$
\sum_{0\le n<N} \delta _n=\inorm{g_N}\le M\mu(g_N,u) = M\mu(u) + 
\half\sum_{0\le n<N} \delta _n ,
$$
whence $\sum_{0\le n<N} \delta _n\le 2M$.  This establishes the desired 
convergence and finishes the proof.

\heading
9. Bump functions and partitions of unity
\endheading

In this section, $\Upsilon $ will be an arbitrary tree. We shall show how to 
define a non-trivial $\C^\infty$ function of bounded support (that is to say, 
a {\sl bump function}) on $\C_0(\Upsilon )$ and how to obtain $\C^\infty$ 
partitions of unity on this space.  We shall use constructions from the 
author's paper [\Haydonfour], from which we quote the following two results.
In the first of them, the set $U(L)$ is the one defined in Section~1 in 
the context of Talagrand operators. 

\proclaim{Proposition 9.1 {\smc (Corollary 3 of [\Haydonfour])}}
Let $X$ be a Banach  space, let $L$ be a set and let $k$ be a positive 
integer or $\infty$. Suppose that there exist 
continuous mappings $S:X\to \ell_\infty(L)$, $T:X\to c_0(L)$ with the 
following properties: 
\roster
\item for all $x\in X$ the pair $(Sx,Tx)$ is in $U(L)\cup \{0\}$;
\item the coordinates of $S$ and of $T$ are all $\C^k$ functions on the sets 
where they are non-zero;
\item $\inorm{Sx}\to \infty$ as $\norm x\to \infty$.
\endroster                      
Then $X$ admits a $\C^k$ bump function.
\endproclaim

\proclaim{Proposition 9.2 {\smc(Theorem 2 of [\Haydonfour])}}
Let $X$ be a Banach space, let $L $ be a set and let $k$ be a positive 
integer or $\infty$. Let $T:X\to c_0(L 
)$ be a function such that each coordinate $x\mapsto T(x)_\gamma $ is of class 
$\C^k$ on the set where it is non-zero. For each finite subset $F$ of $L $,
let $R_F:X\to X$ be of class $\C^k$ and  assume that the following hold:
\roster
\item
for each $F$, the image $R_F[X]$ admits $\C^k$ partitions of unity;
\item
$X$ admits a $\C^k$ bump function;
\item
for each $x\in X$ and each $\epsilon >0$ there exists $\delta >0$ such that
$\norm{x-R_Fx}<\epsilon $ if we set $F =\{\gamma \in L  
:|(Tx)(\gamma )|\ge \delta \}$.
\endroster
Then $X$ admits $\C^k$ partitions of unity.
\endproclaim

\proclaim{Theorem 9.3}
For any tree $\Upsilon $ there is a $\C^\infty$ bump function on 
$\C_0(\Upsilon )$.
\endproclaim
\demo{Proof}
We shall define $L,S,T$ satisfying the conditions of Theorem~9.1. Let $L$ be 
$\Upsilon \times\nats$ and define $S:\C_0(\Upsilon )\to \ell_\infty(L)$ by 
$(Sf)(s,n)=f(s)$. It is obvious that the coordinates of $S$ are $\C^\infty$ 
and that $\inorm{Sf}\to \infty$ as $\inorm f\to \infty$. 
 
In order to define $T$, we first fix a $\C^\infty$ function 
$\phi :\reals^+\to [0,1]$ satisfying
$$
\phi (x)=\left\{\matrix \format \l &\qquad\l \\
                  0 & \text{if $|x|\le \half$}\\
                  1 & \text{if $|x|\ge 1$},
             \endmatrix\right.
$$
and set $\psi =1-\phi $. For $f\in\C_0(\Upsilon )$, $s\in \Upsilon $ and 
$n\in\nats$, we define 
$$
(Tf)(s,n)=\cases
   0\qquad\text{if $f(s)=0$ or there exists $t\in s^+$ with $f(s)=f(t)$}\\
 2^{-n}\phi (2^nf(s))\prod_{t\in s^+}\psi (2^{-n}f(t)/(f(t)-f(s)))
                   \qquad\text{otherwise}.
\endcases
$$

We shall show that the hypotheses of Theorem~9.1 are satisfied.  First we have 
to show that $T$ takes values in $c_0(\Upsilon \times\nats)$. Consider some 
fixed $n$. By Lemma~2.2, for all but finitely 
many $s\in \Upsilon $ there exists $t\in s^+$ with $|f(t)-f(s)|<2^{-2n-2}$. If 
$s,t$ have this property then either $|f(s)|<  2^{-n-1}$, in which case $\phi 
(2^nf(s))=0$, or $|f(s)|\ge 2^{-n-1}$, in which case we have $|f(t)|\ge 2^{-n-
1}-2^{-2n-2}\ge 2^{-n-2}$ and hence $\psi(2^{-n}f(t)/(f(t)-f(s)))=0$.  Thus 
there are, for each $n$, only finitely many $s$ with $(Tf)(s,n)\ne0.$ Combined 
with the inequality $|(Tf)(s,n)|\le 2^{-n}$, this shows that $Tf$  is in 
$c_0(\Upsilon )$. 

Now let $s$ and $n$ be fixed and consider the coordinate mapping $f\mapsto 
(Tf)(s,n)$.  We shall show that every $f_0\in \C_0(\Upsilon )$ has a 
neighbourhood on which this function is $\C^\infty$. If $f_0\in \C_0(\Upsilon 
)$ is such that $f_0(s)=0$ or there exists $t\in s^+$ such that $f(t)=f(s)$, 
then by the calculation in the last paragraph we have $(Tf)(s,n)=0$ whenever 
$\inorm{f-f_0}<2^{-2n-3}$.  So our mapping is certainly infinitely 
differentiable at $f_0$.  Now consider $f_0$ such that $f_0(s)\ne 0$ and 
$\delta =\inf\{|f(t)-f(s)|:t\in s^+\}\ne 0$.  Let $V$ be the set of all $f$ 
with $\inorm{f-f_0}<\sixth\delta$ and let $F$ be the (finite) set of all $t\in 
s^+$ such that $|f_0(t)|\ge\sixth\delta$.  For $f\in V$ and $t\in s^+\less F$ 
we have $|f(t)/(f(t)-f(s))| \le \half$, whence $\psi (2^{-n}f(t)/(f(t)-
f(s)))=1$. Thus, on $V$, our function $(T\cdot)(s,n)$ is the product of 
finitely many $\C^\infty$ functions and is thus itself $\C^\infty$. 

To finish the proof we shall show that for every non-zero $f\in \C_0(\Upsilon 
)$ there is some $(s,n)\in L$ with $(Tf)(s,n)\ne 0$ and 
$|(Sf)(s,n)|=\inorm{Sf}$, that is to say $|f(s)|=\inorm{f}$.  We choose $s$ to 
be maximal in $\Upsilon $ subject to $|f(s)|=\inorm{f}$.  Thus $f(s)\ne0$ and 
there is no $t\in s^+$ with $f(t)=f(s)$.  Let $\delta =\inf\{|f(s)-f(t)|:t\in 
s^+\}$ and choose $n$ so that $2^n\inorm{f}\ge 1$ and $2^{-n}\inorm{f}\le 
\half\delta $.  We then have $\phi (2^nf(s))=1$ and $\psi (2^{-n}f(t)/(f(t)-
f(s)))=1$ for all $t\in s^+$, whence $(Sf)(s,n)=2^{-n}$. \qed\enddemo

\proclaim{Theorem 9.4}
For any tree $\Upsilon $ the space $\C_0(\Upsilon )$ admits $\C^\infty$ 
partitions of unity.
\endproclaim
\demo{Proof}
In 9.2 we take $X$ to be $\C_0(\Upsilon )$ and let $L$ and $T$ be as 
in the proof of Theorem~9.3. We have to define the ``reconstruction 
operators'' $R_F$ and establish the hypotheses (1) and (3) of 9.2. For a 
finite subset $F$ of $L=\Upsilon \times \nats$ we define
$$
(R_Ff)(s) = \left\{\matrix\format \l & \qquad\l\\
          f(s) &\text{if there exists $(t,n)\in F$ with $s\trle t$}\\
           0   &\text{otherwise.} 
               \endmatrix\right.
$$
The image $R_F[X]$ is just the set of continuous functions supported on the 
finite union of intervals $\bigcup_{(t,n)\in F}(0,t]$.  Now this union is 
homeomorphic to some closed interval of ordinals $[0,\Omega]$ and it is shown 
in [\Haydonfour] that the space of continous functions $\C[0,\Omega ]$ admits 
$\C^\infty$ partitions of unity.  Thus hypothesis (1) holds.

Finally, given a non-zero $f\in\C_0(\Upsilon )$ and $\epsilon >0$, we note that 
the set $M$ of maximal elements of $\{t\in\Upsilon: |f(t)|\ge \epsilon\}$ is 
finite.  As in the proof of Theorem~9.3, there exists, for each $t\in M$, a 
natural number $n_t$ such that $(Sf)(t,n_t)=2^{-n_t}$.  If we set $\delta =2^-
{\max_{t\in M}n_t}$ and $F$ is the set of all $(s,m)\in L$ with $(Tf)(s,m)\ge 
\delta $, then $|f(u)|<\epsilon $ for all $u$ in $\Upsilon 
\less\bigcup_{(s,m)\in F}(0,s]$. Hence $\inorm{f-R_Ff}<\epsilon $ as required.
\qed\enddemo

\heading
10. Examples
\endheading

In this section we consider in detail three particular trees, which provide 
most of the counterexamples mentioned in the introduction. Our starting point 
is a  $\reals$-embeddable, non-$\rats$-embeddable tree which is 
well-known to set-theorists. We recall that the standard realization of a full 
countably branching tree of height $\omega _1$ is the set of all functions $t$ 
with domain some countable ordinal and codomain the set 
$\omega $ of natural numbers.  For the rest of this section, $\Lambda $ will 
denote a certain subtree of this tree, namely the set of all {\sl injections} 
$t$ with domain $\dom t$ a countable ordinal and image  $\rg t$ a {\sl 
co-infinite} subset of $\omega$.  This tree is $\reals$-embeddable since the 
function $\lambda :t\mapsto \sum_{\alpha \in\dom t} 2^{-t(\alpha )}$ is 
strictly increasing. We define the $\lambda $-topology $\Cal 
T_\lambda $ on $\Lambda $ by taking basic neighbourhoods of $t$ to be of the 
form $\{u\in [t,\infty): \lambda (u)<\lambda (t)+\epsilon \}$. We start with a 
key lemma about this topology. 

\proclaim{Lemma 10.1}
For the topology $\Cal T_\lambda $ the tree $\Lambda $ is a Baire space.
\endproclaim
\demo{Proof}
It is convenient to work with a slightly different description of the 
$\lambda $-topology.  For $t\in \Lambda $ and $n\in\Omega $ let 
$[t,\infty)_n=\{u\in [t,\infty): t\trl u \text{ and } n\cap \dom u = n\cap 
\dom t\}$. Then we have 
$$
[t,\infty)_{p+1}=\{u\in [t,\infty):\lambda (u)<\lambda (t)+2^{-p}\},
$$
so that the sets $[t,\infty)_n$ form a base of neighbourhoods of $t$ for 
the $\lambda $-topology. 

We shall now show that $(\Lambda ,\Cal T_\lambda )$ is a Baire space by 
showing that it is $\alpha $-favourable [\Choquet]. The strategy for the 
player ($\alpha $) is as follows: if at stage $n$, player $\beta $ plays 
$[t_n,\infty)_{p_n}$ then $\alpha $ chooses $r_n$ to be the $n^{\text{th}}$ 
element of $\omega \less \rg t_n$ and plays $[t_n,\infty)_{q_n}$ where 
$q_n=\max\{p_n,r_n+1\}$. At the end of the game, we have a sequence $t_0\trle 
t_1\trle \cdots$ in $\Lambda $ and there is certainly an ordinal $\alpha 
=\sup_n \dom t_n$ and an injection $t:\alpha \to \omega $ satisfying 
$t\restriction \dom t_n=t_n$ for all $n$. The only question is whether this 
$t$ is in $\Lambda $, that is to say whether $\rg t$ is a co-infinite  subset 
of $\omega $. I claim, in fact, that none of the $r_n$ are in the range of 
$t$; since $\rg t=\bigcap_m \rg t_m$ it is enough to show that $r_n \notin \rg 
t_m$ for all $m$. If $m\le n$ we have $\rg t_m\subseteq t_n$ and $r_n\ne t_n$, 
while for $m\ge n$ $t_m \in [t_n,\infty){q_n}\subseteq [t_n,\infty){r_n+1}$ so 
that $\rg t_m \cap [0,r_n] = \rg t_n \cap [0,r_n]$, which again shows that 
$r_n \notin  \rg t_m$. Finally, since we chose $r_n$ to be the $n^{\text{th}}$ 
member of $\omega \less \rg t_n$, and since $\rg t_n\supseteq \rg t_m$ when 
$m<n$, all the $r_n$ are distinct. \qed \enddemo 

\proclaim{Corollary 10.2}
The $\reals$-embeddable tree $\Lambda $ does not have the Namioka Property.
Hence $\C_0(\Lambda )$ has no Kadec norm, and certainly no LUR norm---though 
by 3.1 it does have a MLUR norm.
\endproclaim
\demo{Proof}
If $t\in \Lambda$ has domain $\alpha $ then the immediate successors of $t$ in 
$\Lambda $ have domain $\alpha +1= \alpha \cup\{\alpha \}$ and are given by 
$t.n$ where 
$$
t.n\restriction \alpha  = t,\quad (t.n)(\alpha )= n.
$$
There is one such successor for each $n\in \omega \less \rg t$. Since $\lambda 
(t.n)\to \lambda (t)$ as $n$ tends to $\infty$ in the infinite set $\omega 
\less \rg t$, we see that $t$ is the $ \Cal T_\lambda $ limit of the sequence 
$t.n$ and, in particular that $\Lambda $ has no $\Cal T_\lambda $-isolated 
points. On the other hand, since $\Cal T_\lambda $ is finer than the reverse 
topology, the map $t\mapsto \indic_{(0,t]}$ is continuous from $(\Lambda ,\Cal 
T_\lambda )$ into $\C_0(\Lambda), \Cal T_{\text p}$. Thus $\Lambda $ does not 
have the Namioka Property as we have just shown that $(\Lambda , \Cal 
T_\lambda )$ is a Baire space. \qed \enddemo

\proclaim{Remark 10.3} {\rm It follows from the above Corollary, together with 
our necessary and sufficient condition for Kadec renormability, that every 
increasing real-valued function of $\Lambda $ has bad points.  There is  a 
direct way to see this.  Let $\rho :\Lambda \to \reals $ be increasing.  For 
any $\alpha \in \reals$ the set $\{t\in \Lambda :\rho (t)>\alpha \}$ is a 
union of wedges $[s,\infty)$, and is hence $\Cal T_\lambda $-open. Thus $\rho 
$ is $\Cal T_\lambda $ lower semicontinuous.  Since $\Lambda, \Cal T_\lambda  
$ is a Baire space, there are points $t$ at which $\rho $ is $\Cal T_\lambda 
$-continuous. Since, for any $t$, the immediate successors of $t$ form a 
sequence which converges to $t$ for $\Cal T_\lambda $, a point at which $\rho 
$ is $\Cal T_\lambda $-continuous must be bad for $\rho $. } \endproclaim 

We shall now use the tree $\Lambda $ considered above to answer negatively the
``three-space problem'' about strictly convex renormings, showing that there 
exists a Banach space $X$ which has a closed subspace $Y$ with a LUR norm,
such that the quotient space $X/Y$ has a strictly convex norm, while $X$ 
itself has no strictly convex renorming.

\proclaim{Proposition 10.4}  
There is a  tree $\Upsilon $ which contains $\Lambda $ as a closed subset, 
with $ \Upsilon \less\Lambda $ discrete, such that $\C_0(\Upsilon )$ admits no 
strictly convex renorming. 
\endproclaim

\proclaim{Corollary 10.5 {\smc(Three-Space Problem)}}
There is a Banach space $X$ and a subspace $Y$ of $X$, isomorphic to 
$c_0(\reals)$, such that the quotient $X/Y$ admits a strictly convex 
renorming, while $X$ does not.
\endproclaim
\demo{Proof of 10.5}
We take $X$ to be $\C_0(\Upsilon )$ and $Y$ to be the subspace consisting of 
all functions $f$ with $f\restriction \Lambda =0$.  The quotient $X/Y$ may be 
identified with $\C_0(\Lambda )$, which admits a strictly convex norm (and 
even an MLUR norm). The subspace $Y$ itself may be identified with 
$\C_0(\Upsilon \less \Lambda) $, which is isomorphic to $c_0(\reals)$ since 
$\Upsilon \less \Lambda $ is discrete and of cardinality the continuum.
\qed \enddemo 

\demo{Proof of 10.4}
For each $t\in \Lambda $ we partition $t^+$ as the union of infinite sets 
$A_1(t)$, $\cup A_2(t)$, and 
augment the tree $\Lambda $ by introducing elements $(t,1)$ and $(t,2)$ with 
the property that $t\trl (t,i)\trl u$ whenever $u\in A_i(t)$. We write 
$\Upsilon $ for the the resulting tree, which contains $\Lambda $ as a closed 
subset. We note that $\Upsilon \less\Lambda = \Lambda \times \{1,2\} $ is open 
and discrete. To show that $\C_0(\Upsilon )$ is not strictly convexifiable we 
shall show that for any increasing function $\phi 
:\Upsilon \to\reals$ there exists $t\in \Lambda $ such that both $(t,1)$ and 
$(t,2)$ are $\phi $-bad, with $\phi (t,1)=\phi (t,2)=\phi (t)$.

Given an increasing function $\phi $ on $\Upsilon $ we define $\psi:\Lambda 
\to\reals$ by $\psi (t)=\phi\restriction\Lambda $. We remarked earlier in 
this section that there exist points at which $\psi$ is  $\Cal T_\lambda 
$-continuous. Let $t$ be any such point, let $i$ be 1 or 2 and let 
$u_n$ $(n\in\omega )$ be a sequence of distinct points of $A_i(t)$. Then we 
have $t\trl(t,i)\trl u_n$, and hence $\psi (t)\le \phi (t,i)\le \psi (u_n)$, 
for all $n$. But the sequence $(u_n)$ converges to $t$ for the topology $\Cal 
T_\lambda $ so that, by continuity of $\psi $ at $t$, we have $\phi (t,i) 
=\phi (t) =\lim_n\phi (u_n)$.  \qed \enddemo 

\smallskip
A very similar construction will now enable us to answer the quotient problem 
for Fr\'echet-differentiable renorming.

\proclaim{Proposition 10.6}Let $\Delta $ be a dyadic tree equipped with an 
increasing function $\rho :\Delta \to \reals^+$ having the property that, for 
all $t\in \Delta $ there is exactly one element $u$ of $t^+$ such that 
$\rho (u)=\rho (t)$. Then $\C_0(\Delta )$ admits a $\C^\infty$ renorming. 
\endproclaim
\demo{Proof}
Since $\Delta $ is finitely branching, the function $\rho $ has no bad points. 
The assumption that each successor set $t^+$ contains only one point $u$ with 
$\rho (u)=\rho (t)$ ensures that $\rho $ is constant on no ever-branching 
subset of $\Delta $.  Thus $\C_0(\Delta )$ admits a $\C^\infty$ renorming by 
Theorem~8.1. \qed\enddemo 

In fact, in the simple case we have just considered, it would have been very 
easy to write down a Talagrand operator and thus establish the proposition 
without recourse to the general result 8.1.  For each $t\in \Delta $ let the 
elements of $t^+$ be denoted $t\dual $ and $\tilde t$, with $\rho (\tilde 
t)>\rho(t)$ and $\rho (t\dual)=\rho (t)$. We set 
$$
(Tf)(t)= (\rho (\tilde t)-\rho (t))(f(t)-f(t\dual)).
$$
It is easy to see that if $t$ is maximal subject to $|f(t)|=\inorm f$ then 
$(Tf)(t)\ne 0$. We show that $\Phi $ takes values in $c_0(\Delta )$ by 
applying the argument used for the operator $S$ in the Lemma~4.3. 
(Indeed, when we note that we here have $F_t=\{\tilde t\}$ for all $t$, we can 
see that the operator $T$ is just $2S$.) 

\proclaim{Proposition 10.7}
The tree $\Lambda$ may be embedded as a closed subset of a dyadic tree $\Delta$
satisfying the hypotheses of Proposition~10.6.
\endproclaim
\demo{Proof}
We augment the tree $\Lambda $ in a way analogous to what we have just done in 
the context of the three-space problem for strictly convex renorming. We 
enumerate, for each $t\in \Lambda $, the elements of $t^+$ (in $\Lambda $) as 
$t_n$ $(n\in\nats)$ and introduce a sequence $t'_n$ of new points in such a 
way that, in the augmented tree $\Delta $, $t$ has exactly two immediate 
successors, $t_0$ and $t'_0$, while for each $n$ the immediate successors of 
$t'_n$ are $t_{n+1}$ and $t'_{n+1}$. We extend the function $\lambda $ to 
$\Delta $ by setting $\lambda (t'_n)=\lambda (t)$ for all $t$ and $n$ and 
observe that the hypotheses of Proposition 10.6 are satisfied. \qed 
\enddemo

\proclaim{Corollary 10.8 {\smc (Quotient Problem for Frechet Renormability)}}
The Banach space $\C_0(\Lambda )$ is a quotient of a space with $\C^\infty$ 
norm, but does not itself admit a Fr\'echet renorming.
\endproclaim
\demo{Proof}
Since $\Lambda $ is a closed subset of $\Delta $,  $\C_0(\Lambda )$ is a 
quotient of $\C_0(\Delta )$, a space which admits a $\C^\infty$ renorming. 
However, $\C_0(\Lambda )$ admits no Fr\'echet renorming, by Theorem~8.1, since 
every increasing real-valued function on $\Lambda $ has bad points. 
\qed\enddemo

\heading 11. Open Problems \endheading

As we remarked in the Introduction, this paper has not offered many results 
about Gateaux differentiability.  Even our sufficient condition for the 
existence of a strictly convex dual norm on $\C_0(\Upsilon )\dual$, established 
in Section~7, does not appear to be a necessary condition.  The only case in 
which the author knows that $\C_0(\Upsilon )$ does not admit a Gateaux smooth 
norm is when $\Upsilon $ is a Baire tree [\Haydonthree].  There is a very big 
gap between this situation and the cases where 7.1 and 8.1 of this paper give 
positive results.  By way of a specific problem, one can ask whether the 
existence of a strictly convex renorming of $\C_0(\Upsilon )$ implies the 
existence of a smooth renorming.  Significant progress in this area will no 
doubt depend on the development of new ways to construct Gateaux-smooth norms.  
To the best of the author's knowledge, the only approaches available at 
present use either strict convexity of the dual norm or Talagrand operators 
(the latter, of course, yielding Fr\'echet-smoothness). 

The other class of problems left open by this paper are those where trees have 
yielded positive results, rather than counterexamples.  In one case it is 
clear that the behaviour of trees is not representative of general locally 
compact spaces: we saw in 5.1 and 8.1 that the existence of an equivalent 
strictly convex or Fr\'echet-smooth norm on $\C_0(\Upsilon )$ implies the 
existence of a bounded linear injection from $\C_0(\Upsilon )$ into a space 
$c_0(\Gamma )$; the Cieselski--Pol space [\DGZ, VII.4.9] shows that $\C(K)$  
may admit a LUR renorming and a $\C^\infty$ renorming but no such injection 
into $c_0(\Gamma )$. The following are problems where more work is needed, 
including perhaps the systematic investigation of new classes of examples. In 
each case, $L$ denotes a locally compact, scattered space.

\proclaim{Problem 11.1}  {\rm Is there any logical connection between LUR 
renormability of $\C_0(L )$ and Fr\'echet-smooth renormability of that 
space?}
\endproclaim

\proclaim{Problem 11.2}  {\rm If $\C_0(L)$ admits a Fr\'echet-smooth 
renorming, does it necessarily admit a $\C^\infty$ renorming?  This question 
seems to be open even in the case where the $\beta ^{\text{th}}$ derived set 
$L^{(\beta )}$ is empty, for some countable ordinal $\beta $ [\DGZ, VII.4].}
\endproclaim

\proclaim{Problem 11.3}  {\rm If $\C_0(L)$ admits a wLUR
renorming, does it necessarily admit a LUR renorming?} \endproclaim

\proclaim{Problem 11.4}  {\rm If $\C_0(L)$ admits a strictly convex
renorming, does it necessarily admit a MLUR renorming?} \endproclaim

\proclaim{Problem 11.5}  {\rm If $\C_0(L)$ is $\sigma $-fragmentable, does it 
necessarily admit a Kadec renorming?} \endproclaim 

Of the above problems, 11.1, 11.3, 11.4 and 11.5 are also open with $\C_0(L)$ 
replaced by a general Asplund space.  To finish, we should perhaps remind the 
reader that the most important open problems in the area of non-separable 
renorming theory are those concerning bump functions and partitions of unity.  
Does every Asplund space admit a $\C^1$ bump function?  Does a space with a 
$\C^1$ bump function necessarily admit $\C^1$ partitions of unity?  In the 
special case of spaces $\C_0(L)$ we may ask whether there always exist 
$\C^\infty$ bump functions.

\enddocument